\documentclass[aip,preprint,11pt]{revtex4-2}

\usepackage{amsfonts, amsbsy, amssymb, amsmath, subfigure}
\usepackage{appendix}
\usepackage[margin = 1in]{geometry}
\usepackage{xcolor}

\definecolor{darkolivegreen}{rgb}{0.33, 0.42, 0.18}

\usepackage{graphicx}
\graphicspath{ {./figures/} }
\usepackage[colorlinks=true,
citecolor=blue,
filecolor=black,
linktoc=all,
linkcolor=magenta,
urlcolor=red]{hyperref}

\usepackage{lineno} 

\begin{document}


\title{Support vector machines for learning reactive islands} 



\author{Shibabrat Naik}
\email[Corresponding author: ]{shibabratnaik@gmail.com}

\author{Vladim{\'i}r Kraj{\v{n}}{\'a}k}
\email[]{v.krajnak@bristol.ac.uk}

\author{Stephen Wiggins}
\email[]{s.wiggins@bristol.ac.uk}
\affiliation{School of Mathematics, University of Bristol \\
Fry building, Woodland Road, Bristol BS8 1UG, United Kingdom}


\date{\today}

\begin{abstract}
	We develop a machine learning framework that can be applied to data sets derived from the trajectories of Hamilton's equations. The goal is to learn the phase space structures that play the governing role for phase space transport relevant to particular applications. Our focus is on learning reactive islands in two degrees-of-freedom Hamiltonian systems. Reactive islands are constructed from the stable and unstable manifolds of unstable periodic orbits and play the role of quantifying transition dynamics. We show that support vector machines (SVM) is an appropriate machine learning framework for this purpose as it provides an approach for finding the boundaries between qualitatively distinct dynamical behaviors, which is  in the spirit of the phase space transport framework. We show how our method allows us to find reactive islands directly in the sense that we do not have  to first  compute unstable periodic orbits and their stable and unstable manifolds. We apply our approach to the H\'enon-Heiles Hamiltonian system, which is a benchmark system in the dynamical systems community. We discuss different sampling and learning approaches and their advantages and disadvantages.
\end{abstract}


\maketitle 


\textbf{Prediction and control of transition between potential wells across index-one saddles are of importance in a diverse array of non-linear systems in physics, chemistry, and engineering. Dynamical systems theory provides the framework for understanding transition dynamics in these systems. The relevant phase space structures are the stable and unstable invariant manifolds of the hyperbolic periodic orbit associated with the index-one saddle. Thus, to determine the fate of an initial condition in the well without solving the non-linear equations, one has to check whether the initial condition lies in a region bounded by the globalized stable (escape from the well) or unstable (entering the well) invariant manifold. The region that bounds the initial condition is given by the intersection of the invariant manifold with a two-dimensional section and called the reactive island, borrowing terminology from chemical reaction dynamics. In this paper, we develop and verify a trajectory-based framework using support vector machines, by applying it to the H{\'e}non-Heiles Hamiltonian, for learning the reactive islands. Our results show that support vector machines are an ideal data-driven framework for learning the geometry of phase space structures. The approaches developed here are robust to changes in system parameters and geometry of the reactive islands.}

\section{Introduction}
\label{intro}

The goal of this paper is to develop a machine learning framework that identifies the phase space structures governing phase space transport in data sets constructed from trajectories of Hamilton’s equations.  Our focus in this paper will be on two degrees-of-freedom Hamiltonian systems.

The Hamiltonian function typically has the form of the sum of a kinetic energy, a function of the momentum variables and the potential energy, a function of the position variables. For Hamilton’s equations in canonical form each momentum variable is canonically conjugate to one position variable. Hence there are an equal number of momentum and position variables. The space of both the momentum and position variables is referred to as the phase space and the space of only the position variables is referred to as the configuration space.  The number of degrees of freedom of a Hamiltonian system is the number of configuration space variables. 

Hamilton's equations describe dynamics in phase space. Nevertheless, the topography of the potential energy surface plays an important role in defining the reaction dynamics problem for a specific Hamiltonian system. In particular, potential wells are identified as reactants and products. Wells are separated by saddle points and the configuration space picture of reaction dynamics is that of a trajectory evolving between wells by ``crossing'' the saddle point.

Index-one saddle points on a potential energy surface (PES) are the ``seeds'' for the phase space structures from which the theory of reactive islands is constructed. Conley \cite{conley1969ultimate, conley1968low} was the first to analyze the phase space geometry and dynamics near an index-one saddle in two degrees-of-freedom (DoF) Hamiltonian systems in his studies of the circular restricted three body problem.

The Lyapunov subcenter theorem \cite{moser1958generalization, kelley1967liapounov} is fundamental for passing from an equilibrium point (the index-one saddle) to dynamical behavior. It states that for a range of energies above the energy  of the index-one saddle (the exact range is not given in the theorem)  there exists an hyperbolic periodic orbit having two-dimensional stable and unstable manifolds. In the three-dimensional energy surface, stable and unstable manifolds have the geometry of cylinders (``tubes'') and these stable and unstable cylinders mediate transport across the region of the index-one saddle.

The cylindrical structure of the stable and unstable manifolds, together with their invariance, implies that trajectories starting within the tubes must remain in the tubes for all positive and negative time. All trajectories inside the stable cylinder approach the hyperbolic periodic orbit in positive time, pass through the region bounded by the periodic orbit, and then exit the region through the unstable cylinder. The stable cylinder (and the unstable cylinder) has two ``branches'' joined at the periodic orbit, one emanating to each ``side'' of the orbit.

In the context of escape from a potential well, escaping trajectories must lie in a branch of the stable cylinder. Trajectories starting in the branch of the stable cylinder lying in the potential well correspond to forward escape trajectories, i.e. trajectories that escape in forward time. Trajectories starting in the branch of the stable cylinder lying outside the potential well correspond to backward escape trajectories, i.e. trajectories that are captured in the well in forward time. After escape from/capture in the well, trajectories are guided by an unstable cylinder. This is how stable and unstable cylinders govern the dynamics in phase space. 

We describe a lower dimensional technique for probing the geometry of the cylinders. In the potential well we construct a two dimensional Poincar\'e section, i.e. a two dimensional surface where the Hamiltonian vector field is everywhere transverse to the surface. The Poincar\'e section is constructed in such a way that the stable (and unstable) cylinder intersects it in a topological circle. The region bounded by this topological circle is referred to as a ``reactive island'', where `reactive' refers to the occurence of a chemical reaction, an analogue of the escape from a potential well considered here. The Poincar\'e map of the Poincar\'e section into itself is the map that associates to a point its first return to the Poincar\'e section under the flow generated by the Hamiltonian vector field. The inverse of this Poincar\'e map of the Poincar\'e section into itself is the map that associates to a point its first return to the Poincar\'e section under negative time, i.e. the point ``where it came from''.

We consider the preimages of this reactive island by considering its evolution backwards in time under the inverse of the Poincar\'e map. In this way one obtains a reactive island on the Poincar\'e section that returns to the original reactive island, and then escapes for positive time evolution. By repeating this construction we obtain a sequence of reactive islands, ordered in time,  which sequentially map to each other before reacting. Possible  pathological  intersections of the cylinders with the Poincar\'e section can occur, and are discussed in \cite{deleon1989}. A theory of reaction dynamics for two DoF systems based on the geometry of these stable and unstable cylinders was developed in the late 1980's and 90's in \cite{marston1989,deleon1989,OzoriodeAlmeida90,deleon91a, deleon91b, deleon92, mehta1992classical, fair1995}, although some of the ideas appeared in earlier work \cite{DeVogelaere1955,PollakChild80}, and it goes by the name of ``reactive island theory''.  

Our goal is to consider a data set consisting of points on the energy surface that are  labeled based on the evolution of the corresponding trajectory
and ``learn'' which points lead to a escape from a potential well and which do not.  We have chosen to use support vector machines (SVM), a class of supervised learning classification and regression algorithms \cite{Cortes95, Vapnik96,Vapnik2000}, which has been introduced to nonlinear dynamical systems by Ref. \onlinecite{Mukherjee1997}. We employ the classification algorithms, which define a decision function by determining the boundary between different classes of data. To be precise, SVM identifies a subset of the data set referred to as ``support vectors'' to calculate the boundary for which the distance to data in both classes is locally maximal. In our case, the classes in the classification correspond to initial conditions leading to ``qualitatively different'' dynamical behaviour. The learned boundary between escaping and non-escaping trajectories in phase space consists of the invariant manifolds discussed above and thereby SVM enables us to determine the geometrical structures in phase space governing reaction dynamics.

The reasons for choosing support vector machine (SVM) to identify the phase space structures are: (i) Nonlinear kernels in SVM provide means to approximate curves such as reactive islands which form nonlinear boundaries between reactive and non-reactive trajectories. We would like to note that methods for nonlinear clustering~\cite{wang_nonlinear_2016} and other kernel methods~\cite{scholkopf_statistical_2000} are also candidates for developing similar approaches, but beyond the scope of this study. (ii) We design our approach with extensions to higher-dimensional applications in mind. SVM is known to work well even with small amounts of data, therefore our approach is computationally better suited than existing methods for high-dimensional systems and systems where integrating trajectories is expensive. This makes SVM suitable for the reactive island theory of three DoF~\cite{krajnak2021reactive} and system-bath models of isomerization~\cite{naik2020} which has been developed recently and supports the case for a computationally efficient approach to phase space structures presented here.

The Hamiltonian system that we will use to benchmark our approach will be the H\'enon-Heiles system \cite{henon_applicability_1964}. This is a two degree of freedom Hamiltonian system that serves as a paradigm for understanding complex dynamics in a variety of settings. As a function of the energy, it can display dynamical behavior that numerically appears to be ``near integrable'' to completely chaotic. Moreover, this system has three index one saddles that define three distinct reaction channels having the geometric structure discussed above. The geometry of reaction dynamics geometry for a similar system has been analyzed and discussed in Ref.~\cite{naik2019b}.

This paper is outlined as follows. In Section \ref{sec:HH} we describe the two degrees-of-freedom  H\'enon-Heiles Hamiltonian system and the nature of reaction dynamics in the context of this model. In Section \ref{sec:SVM} we describe the machine learning technique of support vector machines (SVM) and the approach of active learning. We describe how trajectories are constitute into data sets and the use of Lagrangian descriptors as a trajectory diagnostic. In Section \ref{sec:results} we describe our results, and in Section \ref{sec:concl} we present our conclusions and the outlook for further research.

\section{H{\'e}non-Heiles system and reaction dynamics}
\label{sec:HH}

We use the H{\'e}non-Heiles Hamiltonian with three index-one saddles in bottlenecks through which trajectories can escape. This escape can be interpreted as a reaction by crossing the potential energy barrier. This model system and its high dimensional analog has been studied in great detail in nonlinear dynamical systems, statistical mechanics, for developing molecular simulation algorithms~\cite{noid_spectral_1977,pullen_comparison_1981,berdichevsky_statistical_1991,contopoulos_systems_2005,zhao_threshold_2007}. In this study, we will define the three exits as follows: the entering via the top index-one saddle corresponds to the formation of a molecule complex by combination of two atoms or molecules, while the left and right exits correspond to dissociation of the molecule into two different products with structural symmetry.

The Hamiltonian is given by
\begin{align}
	\mathcal{H}(x, y, p_x, p_y) & =  T(p_x, p_y) + V_{\rm HH}(x, y) \nonumber \\
	& = \frac{1}{2m_x}p_x^2 + \frac{1}{2m_y}p_y^2 + \frac{1}{2}\omega_x^2 x^2 + \frac{1}{2}\omega_y^2 y^2 + x^2 y - \dfrac{\delta}{3}y^3
	\label{eqn:henon-heiles}
\end{align}
where all the parameters are set to $1.0$ to be comparable with the known results in the literature~\cite{henon_applicability_1964,henon_integrals_1974,demian2017}. The vector field is given by

\begin{equation}
	\begin{aligned}
		\dot{x} = & \dfrac{p_x}{m_x} \\
		\dot{y} = & \dfrac{p_y}{m_y} \\
		\dot{p}_x = & - \omega_x^2 x - 2 x y \\
		\dot{p}_y = & - \omega_y^2 y - x^2 + \delta y^2 
	\end{aligned} 
	\label{eqn:hh2dof_vf}
\end{equation}

and the equilibrium points are

\begin{align}
	\left( 0, 0, 0, 0 \right), \qquad \left(0, \dfrac{\omega_y^2}{\delta}, 0, 0 \right), \qquad \left( \pm \omega_x \sqrt{\dfrac{\omega_x^2}{2} + \dfrac{\delta \omega_x^2}{4}}, - \dfrac{\omega_x^2}{2}, 0, 0 \right).
\end{align}

The eigenvalues of the linearized system in Appendix~\ref{eqn:linearized_hh2dof} at the equilibrium points gives the origin is a center $\times$ center equilibrium and the three remaining are saddle $\times$ center equilibrium points.

\begin{figure}[!ht]
	\centering
	\includegraphics[width=0.45\textwidth]{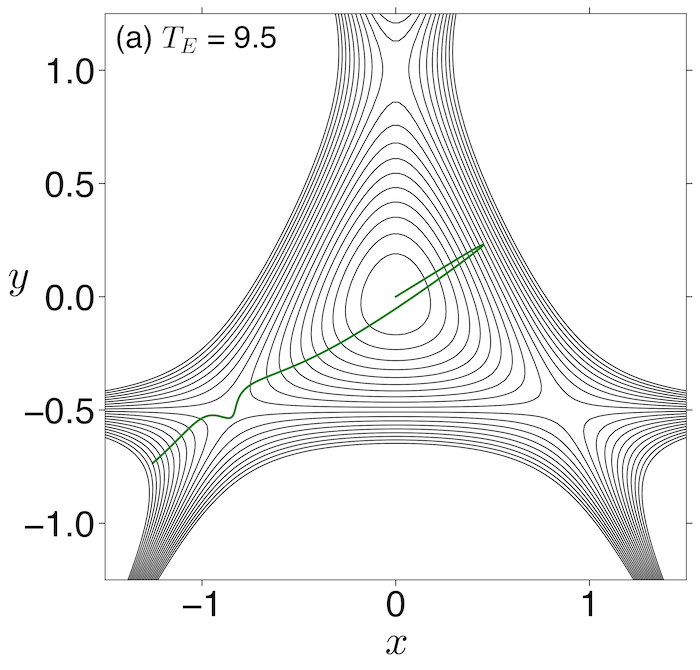} 
	\includegraphics[width=0.45\textwidth]{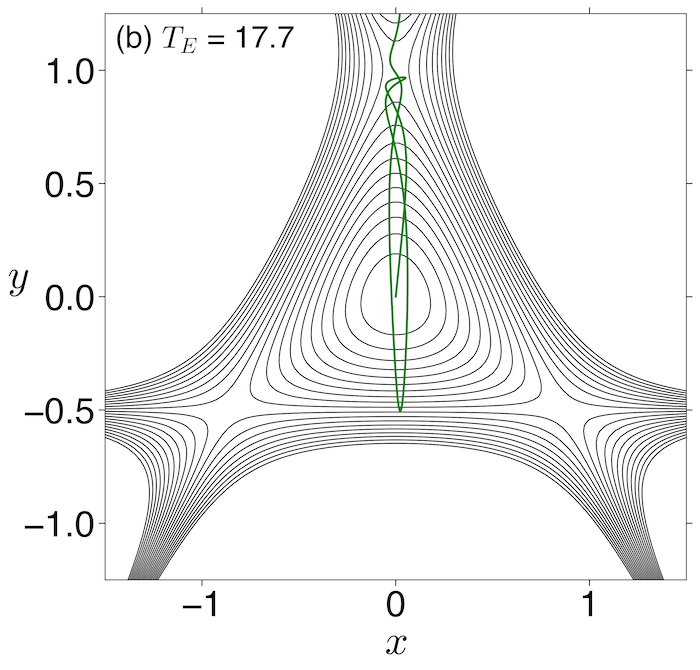}
	\includegraphics[width=0.45\textwidth]{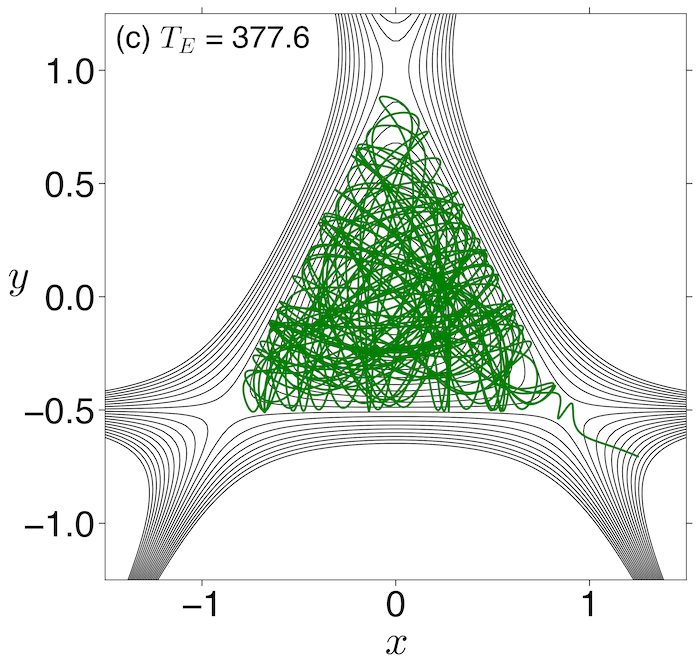}
	\includegraphics[width=0.45\textwidth]{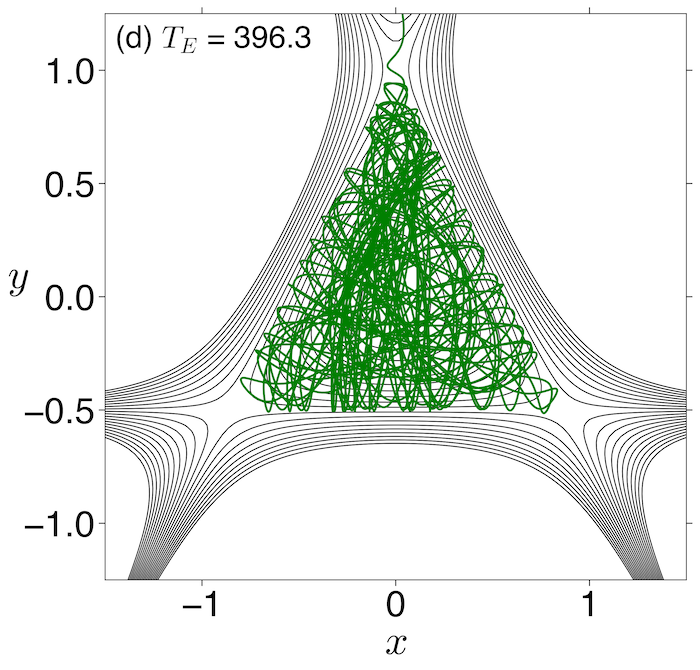}
	\caption{Four sample trajectories initialized on the section $y = 0, p_y > 0$ with $p_x = 0.516,0.07,0.526,0.08$ in (a-d), respectively, projected on the configuration space. The total energy $E = 0.17$ is slightly above the energy of the index one saddles and escape times $T_E$ are shown on each plot. }
	\label{fig:pes_proj_HH}
\end{figure}

We illustrate the dynamics given by the vector field~\eqref{eqn:hh2dof_vf} by showing four sample trajectories in the Fig.~\ref{fig:pes_proj_HH} at same energy starting at the center $\times$ center equilibrium point with different $p_x$ momenta and $p_y > 0$. In Fig.~\ref{fig:pes_proj_HH} (a,c), $p_x$ coordinates only differ in the second significant digit, and yet the escape from the potential well occurs via different bottlenecks and at different times. However, in Fig.~\ref{fig:pes_proj_HH} (b,d), a similar difference in $p_x$ coordinates only changes the escape time. This illustrates the challenge of sampling an ensemble of escape (transition) trajectories. Furthermore, the challenge in identifying different timescales of escape trajectories is due to the large variation in time to escape from the potential well when there is merely a small difference in $p_x$. We will next show the underlying phase space structure of such escape behavior and then show the use of trajectory data in identifying the structure. Finally, it is important to note that for the total energies that we consider the energy surface is unbounded. This implies that notions of recurrence from standard ergodic theory, such as Poincar\'e recurrence and ergodicity, do not apply because escaping trajectories become unbounded and never return.

\begin{figure}[!ht]
	\centering
	(a)\includegraphics[width=0.49\textwidth]{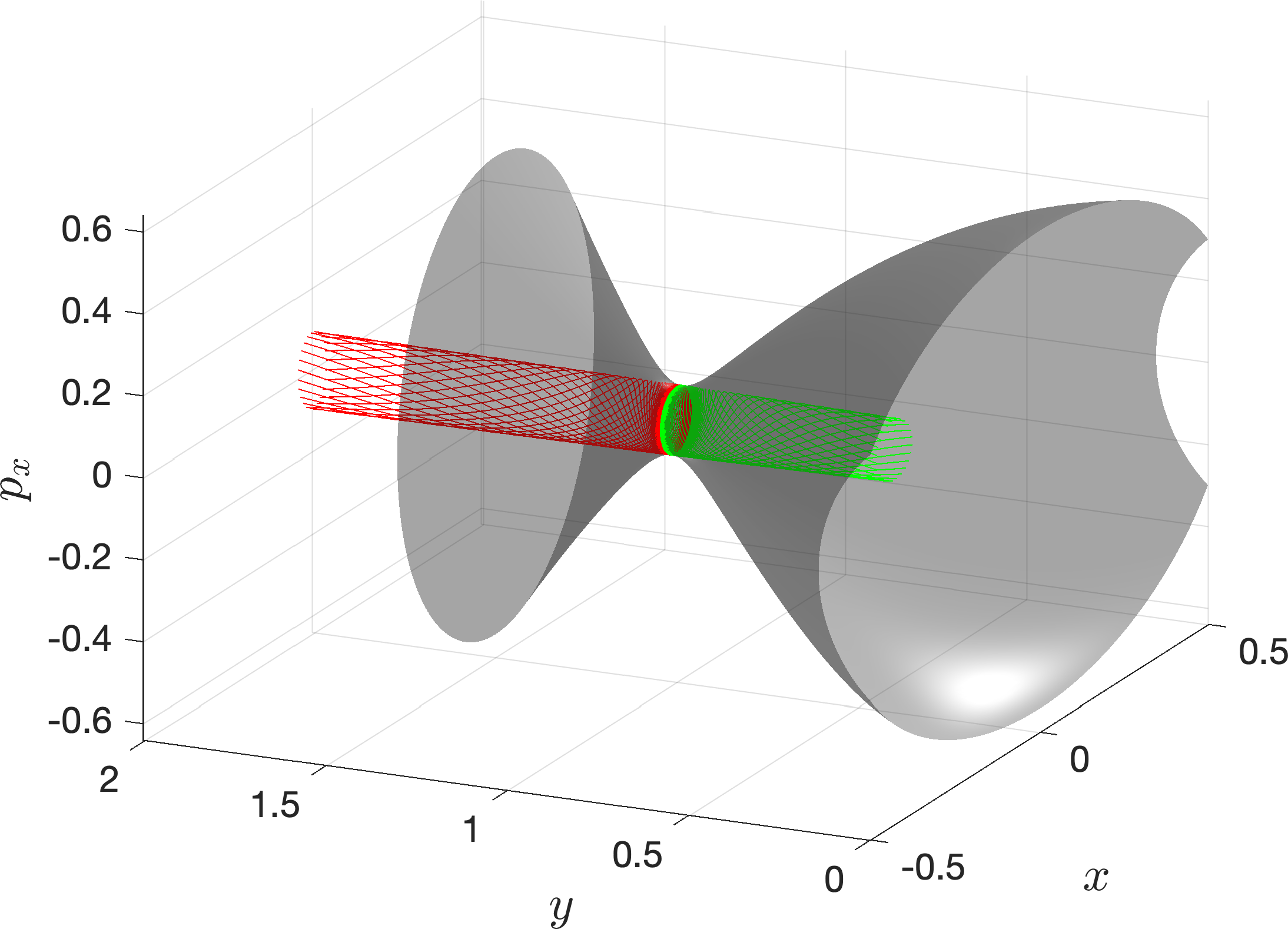}
	(b)\includegraphics[width=0.42\textwidth]{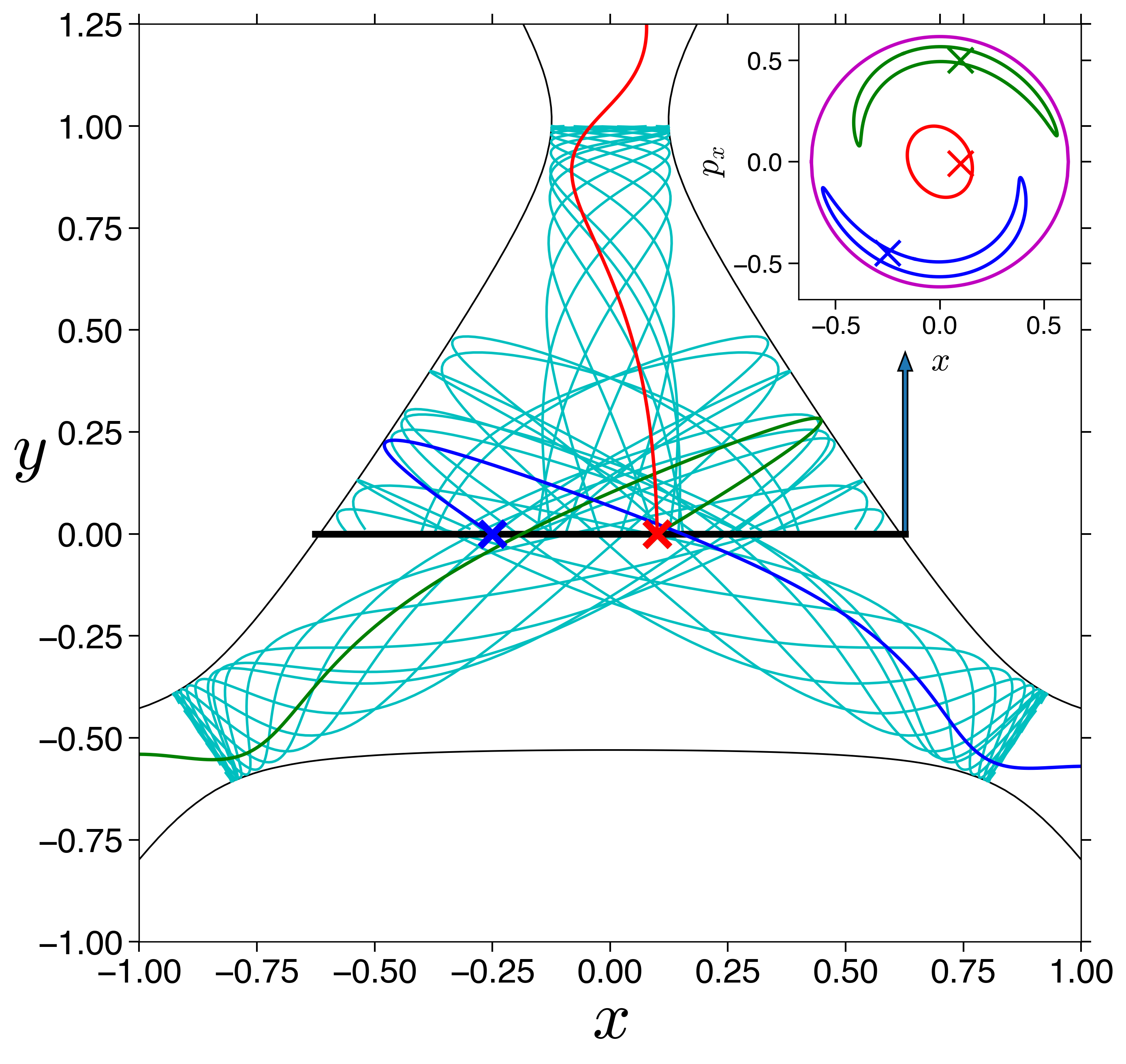}
	\caption{ (a) Cylindrical (or tube) manifolds, stable in green and unstable in red, of the hyperbolic periodic orbits associated with the top index-one saddle in the H{\'e}non-Heiles Hamiltonian. The energy of the hyperbolic periodic orbit and the invariant manifolds are at the total energy, $E = 0.17$ and mediate the trajectories that escape via top saddle as shown in Fig.~\ref{fig:pes_proj_HH}(b,d). (b) Stable manifolds projected on the configuration space reveal the geometry of imminent escape from the potential well via the three bottlenecks. Only the segment of the stable manifolds from the hyperbolic periodic orbits to the intersection with the Poincar{\'e} section (shown as a black line) is shown for the energy, $E = 0.19$.}
	\label{fig:manifolds_HH}
\end{figure}

The trajectory behavior shown in Fig.~\ref{fig:pes_proj_HH} is mediated by the stable manifolds (cylindrical geometry or tubes) of the hyperbolic periodic orbits associated with the index-one saddles. For the two degrees-of-freedom Hamiltonian system, tube manifolds can be computed and visualized in their complete geometry as shown in Fig.~\ref{fig:manifolds_HH} (a) for the index-one saddle at $\left(0, \dfrac{\omega_y^2}{\delta}, 0, 0 \right)$. The stable manifolds associated with the three saddles are projected on the configuration space $(x,y)$ are shown in Fig.~\ref{fig:manifolds_HH} (b) (In Appendix.~\ref{app:tube_manifolds}, Fig.~\ref{fig:manifolds_all_saddles} shows all the tube manifolds in 3D). In this article, we compute the stable manifolds using the procedure described in the Appendix~\ref{app:tube_manifolds} to obtain the segments which intersect a section with $y = 0$ and $p_y > 0$ shown in Fig.~\ref{fig:manifolds_HH}(b). Given such a section of the three dimensional energy surface, the first order reactive islands of escape are defined as the last intersection of the stable manifolds with the section. By last intersection, we mean the trajectories in forward time intersect the section and then leave the potential well without returning to the section and also referred to as the {\it imminent escape}.

\section{Support vector machines and active learning}
\label{sec:SVM}

The classification algorithms for SVM \cite{Cortes95,Vapnik96,Vapnik2000} construct a boundary between different classes of data. In our case, the classes correspond to areas of qualitatively different dynamics on the section $y = y_c$ and $p_y > 0$. By qualitatively different dynamics we mean either areas leading to imminent escape from the potential well over the corresponding saddle or the complementary area that does not lead to imminent escape. The exact boundaries between these areas are the reactive islands formed by stable and unstable invariant manifolds of hyperbolic periodic orbits discussed above. A SVM classification algorithm, also referred to as support vector classifier (SVC), therefore approximates reactive islands in this setting.

\begin{figure}[!h]
	\centering
	\includegraphics[width=0.75\textwidth]{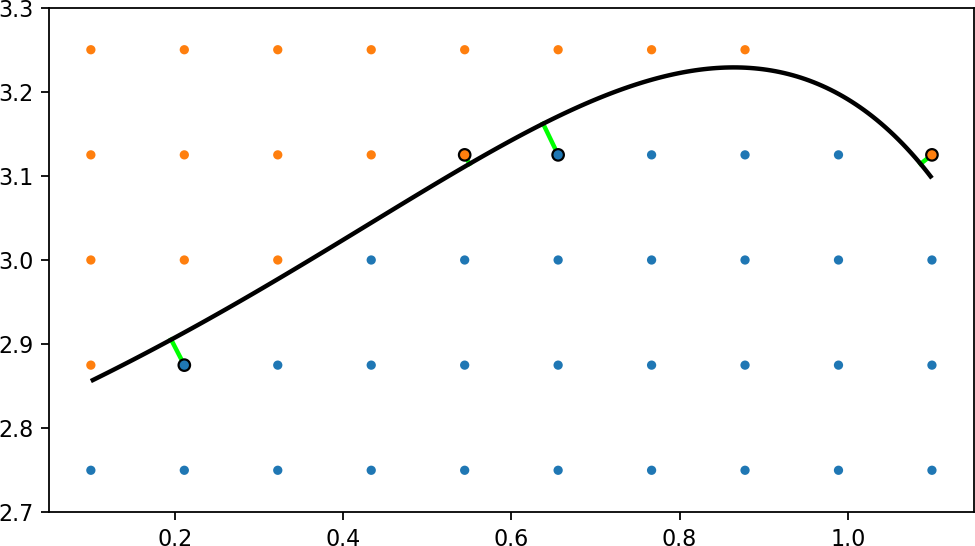}
	\caption{An illustration of a decision boundary (black) between two classes of data (blue and orange) calculated by SVC with radial basis function kernel. The distance between the boundary and the closest points $P_i$ of every class, in this case the support vectors, is highlighted in green.}
	\label{fig:svm_illustration}
\end{figure}

Similarly to~\onlinecite{pozun2012optimizing}, we use the \verb|scikit-learn| \cite{scikit-learn} implementation of SVM \cite{Chang2011}. The implementation can be used with various kernels, of which the radial basis function kernel is best suited to approximate reactive islands in the H\'enon-Heiles system, which are topological circles. With this kernel, a previously unseen data point $P$ is predicted to belong to a class using the decision function
\begin{equation}
	\sum\limits_i \alpha_i l_i e^{-\gamma ||P_i-P||^2},
\end{equation}
where $\gamma>0$ controls the width of the Gaussian, $l_i=\pm1$ are class labels of training data $P_i$ and $C\geq\alpha_i\geq0$ are weights calculated by the algorithm, of which only a number of the weights $\alpha_i$ is non-zero. These weights correspond to the training data $P_i$ called support vectors. The weights $\alpha_i$ are calculated by SVC such that the distance between the predicted boundary and the closest points $P_i$ of every class is maximised, as illustrated in Fig. \ref{fig:svm_illustration}.

\begin{figure}[!h]
	\centering
	\includegraphics[width=0.75\textwidth]{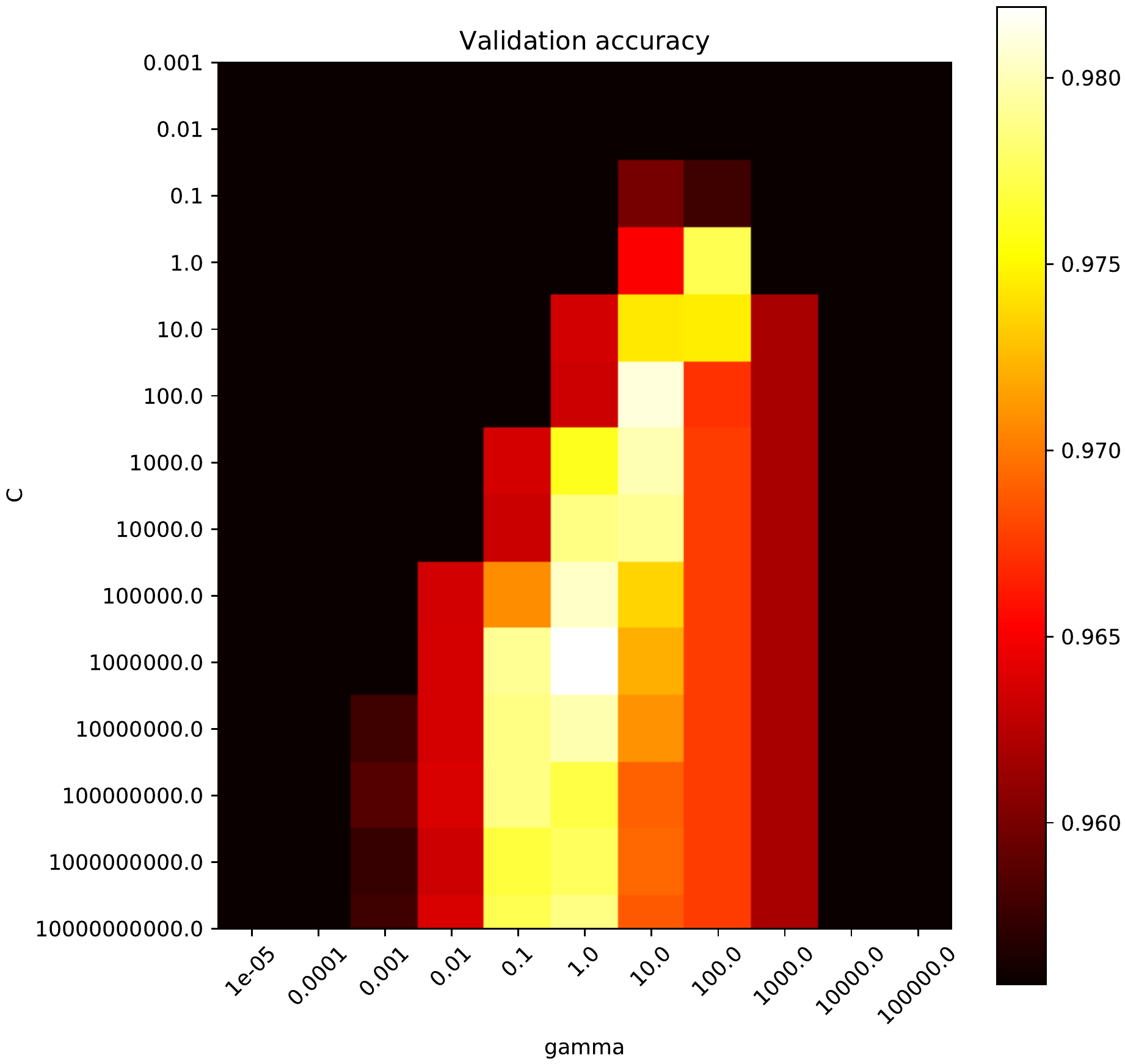}
	\caption{Colormap showing accuracy for different combination of radial basis function parameters, $(C,\gamma)$. The accuracy is obtained using a 5-fold cross validation over the grid of values for $C$ and $\gamma$. The pair of value which gives maximum accuracy is chosen for training the support vector classifier.}
	\label{fig:henonheiles2dofrbfparametersoptimizationset2}
\end{figure}

The upper bound $C$ on weights $\alpha_i$ is a user defined value that controls the complexity of the decision boundary - a low value of $C$ gives a smoother decision boundary, while a high value of $C$ leads to higher accuracy. In this article, we first perform a search over a wide interval of $C,\gamma$ values as shown in Fig.~\ref{fig:henonheiles2dofrbfparametersoptimizationset2}. Then, a smaller interval for both parameters is chosen for each of the support vector classifier approaches. The cross-validation ensures that the trained model does not suffer from over-fitting by splitting the training data into 5 folds, each of which is used as a test set with the remaining four as training set.

While it is possible to apply SVM to a fixed training data set, such as a regular grid, the accuracy of the resulting decision boundary will be limited by the amount of training data and its spacing. A significantly less data-intensive approach is offered by active learning \cite{settles09,Kremer14}, where the `learner' biases its sampling based on information obtained from previous samples. To do this, we start SVM with a coarse grid of data points and iteratively add data points in the proximity of the decision boundary and re-run SVM. This allows the algorithm to explore the intricate structures usually formed by invariant manifolds in systems describing chemical reaction dynamics. At every iteration we randomly add one data point near $10$ randomly selected support vectors. The point is added using the multivariate normal distribution $\mathcal{N}(P,I)$, where $P$ is the support vector and $I$ the identity matrix.

We would like to point out the importance of the precise problem formulation. Homoclinic and heteroclinic intersections of invariant manifolds lead to fractal structures, that is a fractal boundary between classes of dynamics. There is no known way to resolve fractal structures with finite precision and a finite number of data points. Thus approximating the boundary using fixed-width Gaussians is bound to fail. In many systems it is possible to avoid these fractal structures by carefully selecting a surface of initial conditions and studying the dynamics under the corresponding return (Poincar\'e) map.

\section{Results and discussion}
\label{sec:results}

In the figures below, the cyan dots denote the support vectors used by the classifier to learn the decision boundary and the black dashed line denotes the learned decision boundary. The true reactive islands are shown as red, green, and blue curves.

\textbf{Fixed training data.} We sample initial conditions on the two-dimensional section of the three dimensional energy surface, $\mathcal{H}(x, y, p_x, p_y) = E$, with $y = y_c$ and $p_y > 0$. In this study, we present the results for two sections at $y_c = 0, -0.25$ to denote distinct sections at the location of the well and near the bottleneck. We also show the results for $E = 0.17-0.20$ in increments of $0.1$ to illustrate the approach for various imbalance (ratio of reactive to non-reactive trajectories) in the training data corresponding to different excess energies.   
We generate a grid of initial conditions $100 \times 100$ and sample the $y-$momenta using the fixed energy condition. Then, we run trajectories with the initial conditions for a prediction time horizon of $t = 30$ time units. Then, we classify the escape trajectories as reactive through bottleneck 1 or 2 or 3 if they cross the line $x = -1.25$, or $x = 1.25$, or $y = 1.25$, respectively. If an initial condition does not satisfy any of the above conditions for the chosen time interval, we label it as non-reactive denoted by 0. Thus, we obtain a multi-label training dataset with two feature vectors, $x,p_x$ coordinates to discover the distribution of reactive trajectories on the two-dimensional section corresponding to the two coordinates. We use a smaller interval for the Gaussian radial basis function parameters, $C \in \{1e2, 1e3, 1e4, 1e5\}$ and $\gamma \in \{10, 1e2, 1e3, 1e4\}$, to perform a grid search for high accuracy along with 5-fold cross-validation during the training. 

\begin{figure}[!h]
	\centering
	\subfigure[$E = 0.17, y_c = 0$]{\includegraphics[width = 0.24\linewidth]{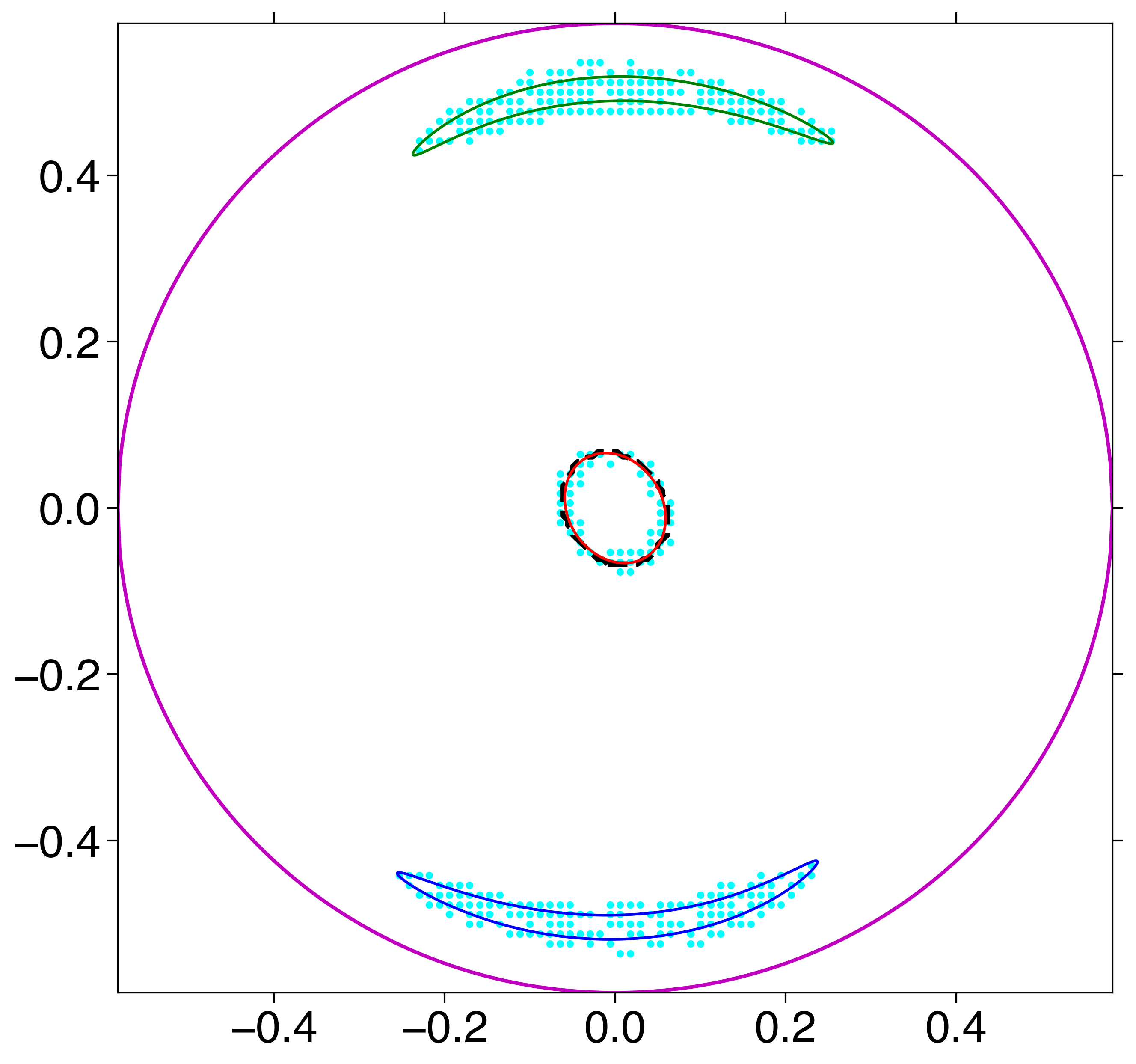}}
	\subfigure[$E = 0.18, y_c = 0$]{\includegraphics[width = 0.24\linewidth]{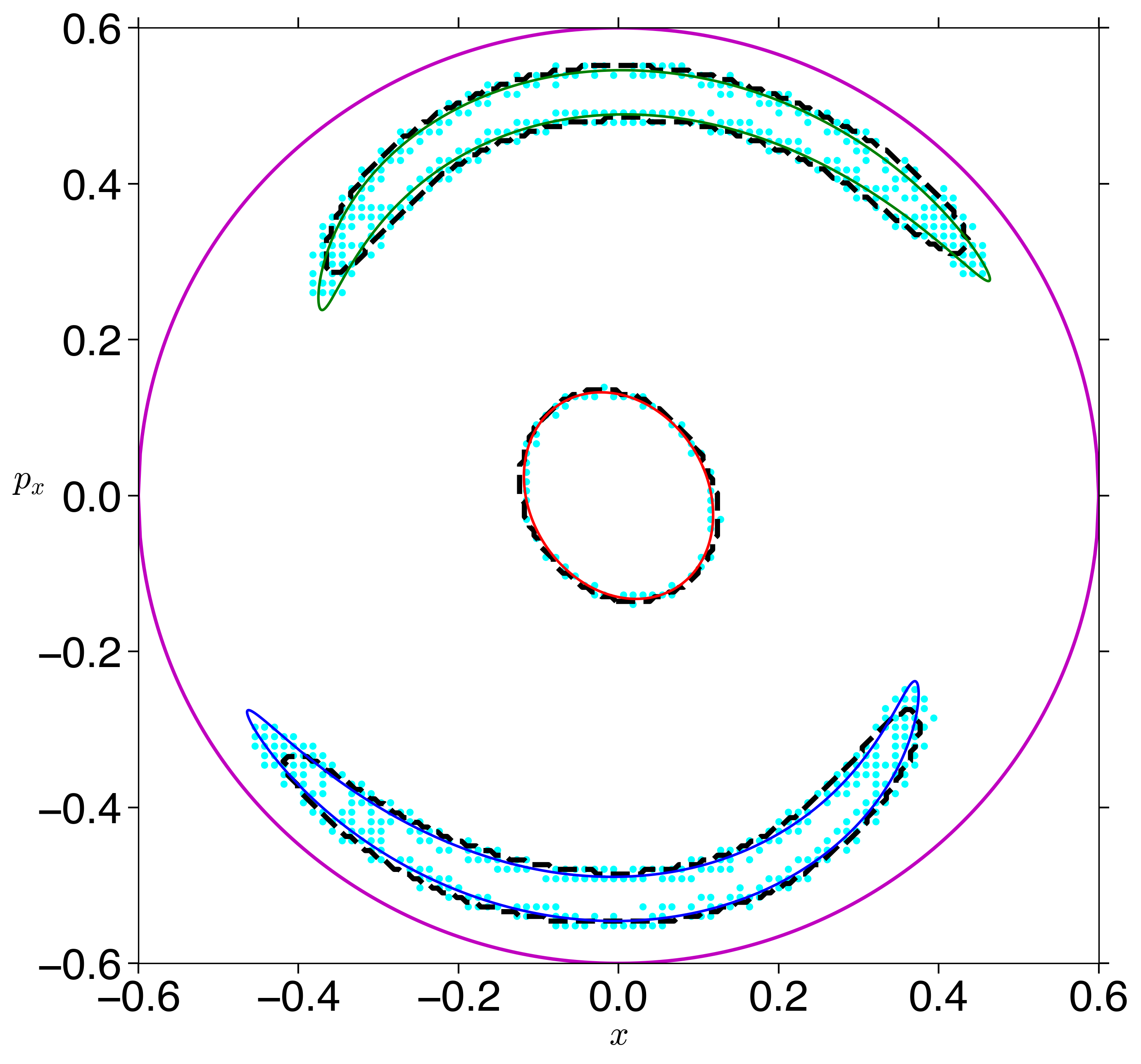}}
	\subfigure[$E = 0.19, y_c = 0$]{\includegraphics[width = 0.24\linewidth]{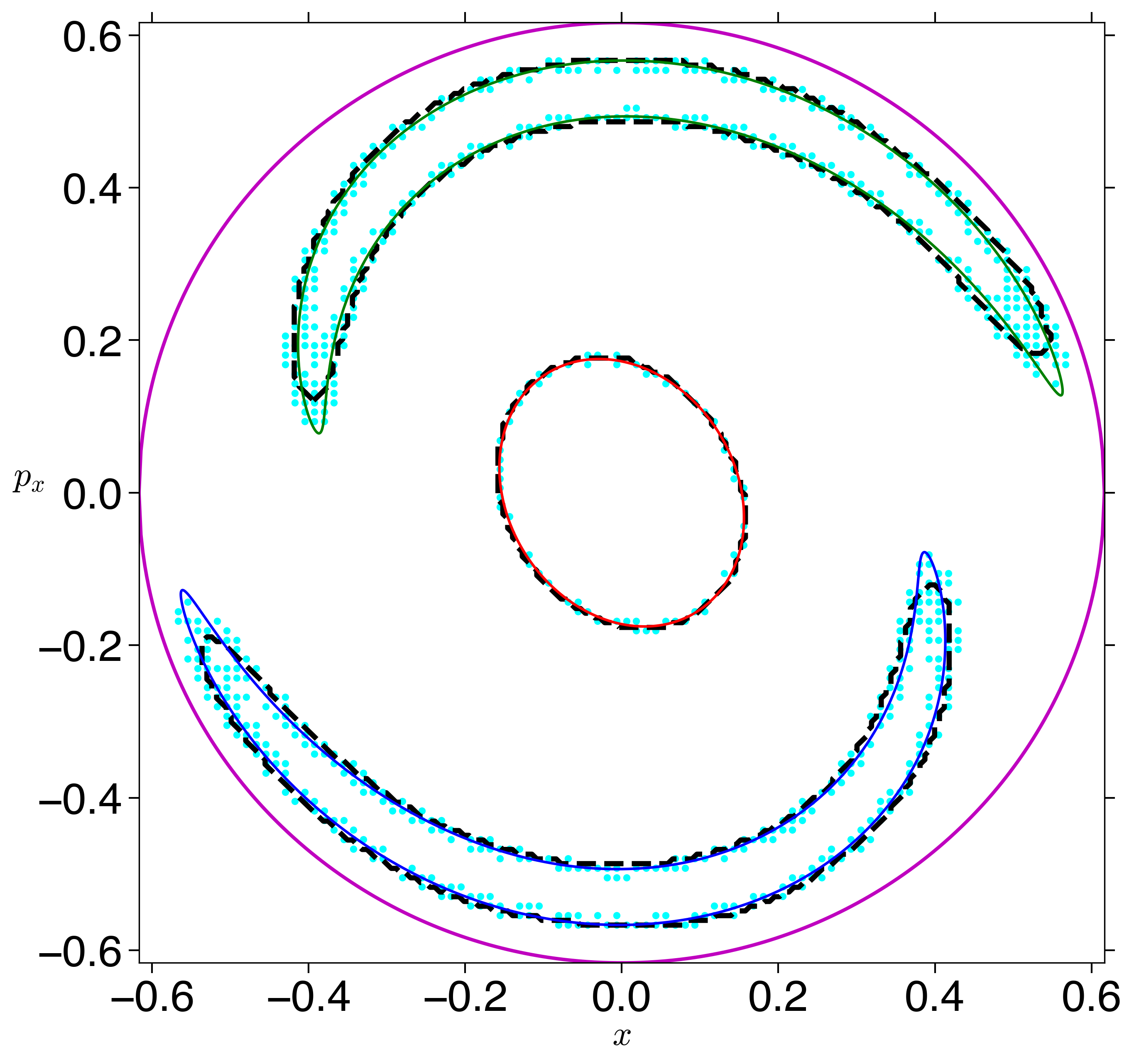}}
	\subfigure[$E = 0.20, y_c = 0$]{\includegraphics[width = 0.24\linewidth]{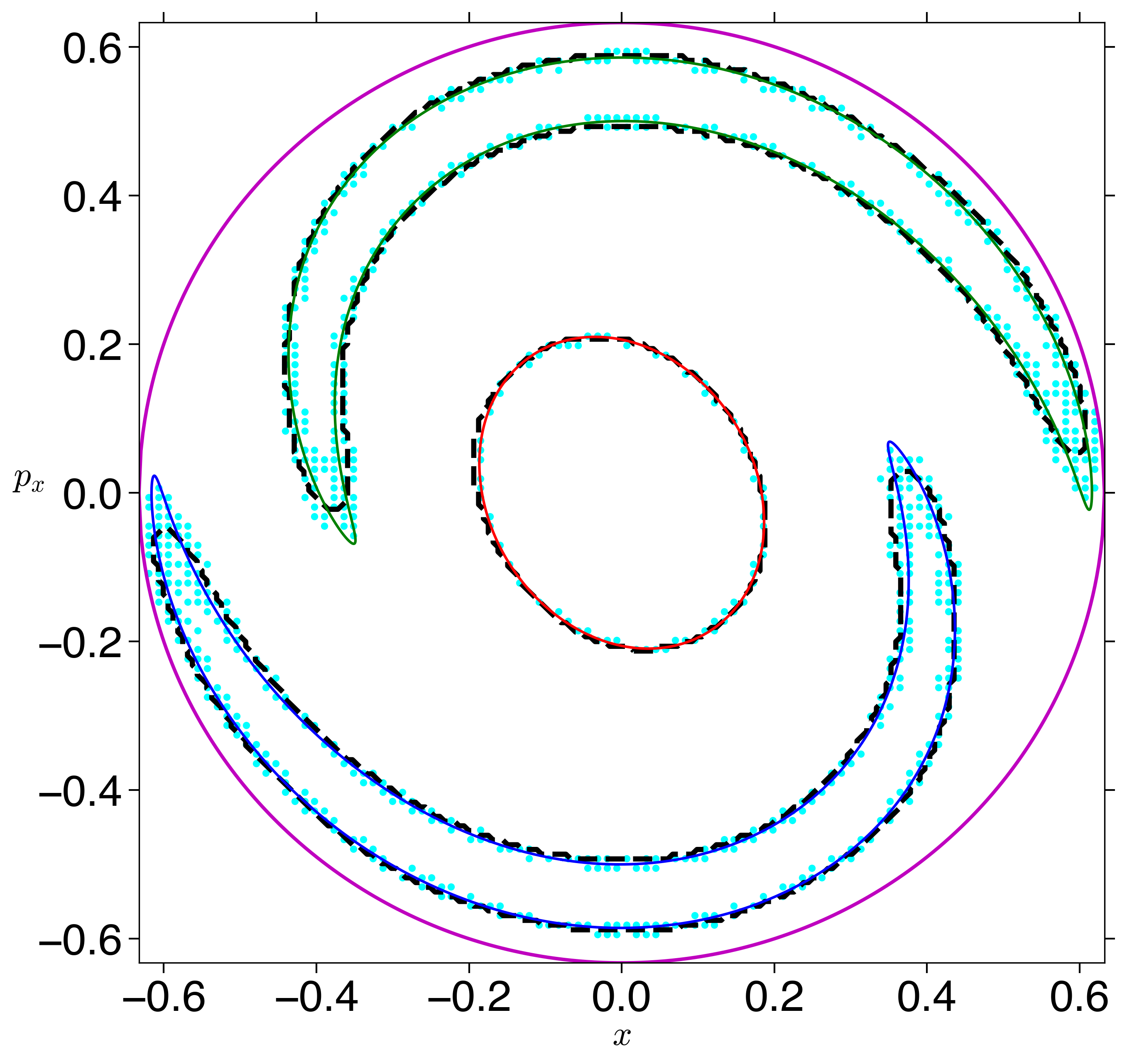}}
	\subfigure[$E = 0.17, y_c = -0.25$]{\includegraphics[width = 0.24\linewidth]{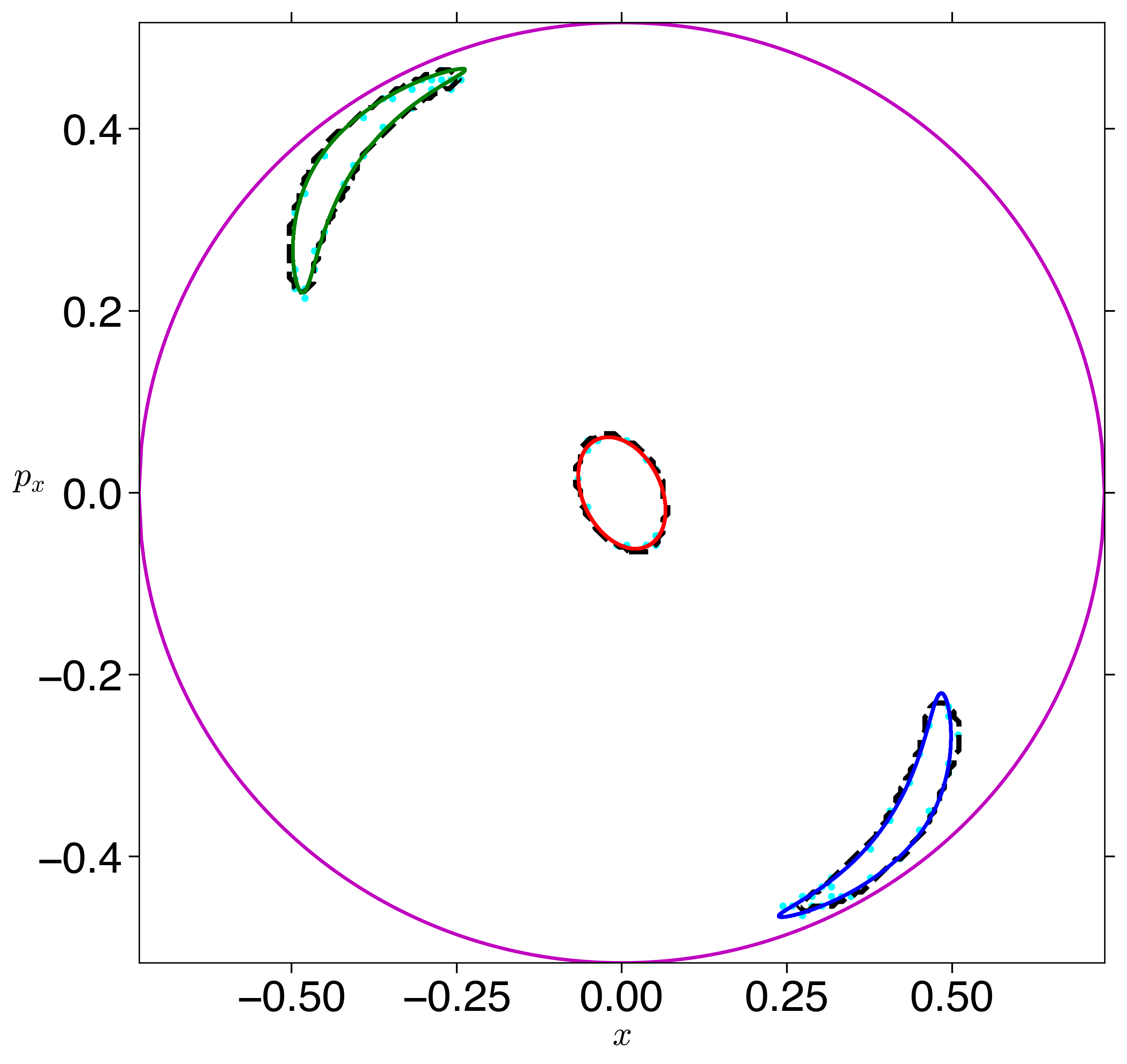}}
	\subfigure[$E = 0.18, y_c = -0.25$]{\includegraphics[width = 0.24\linewidth]{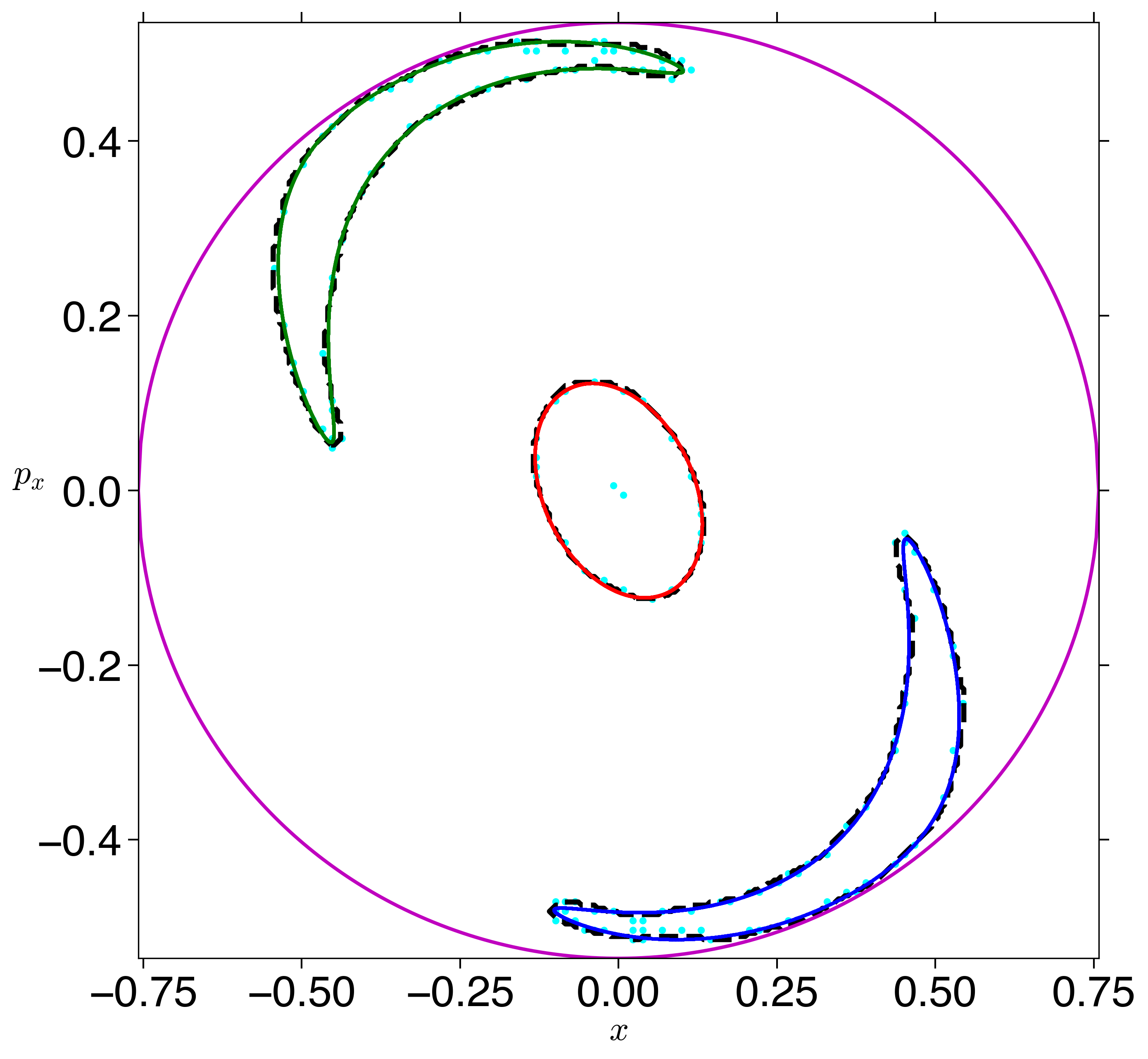}}
	\subfigure[$E = 0.19, y_c = -0.25$]{\includegraphics[width = 0.24\linewidth]{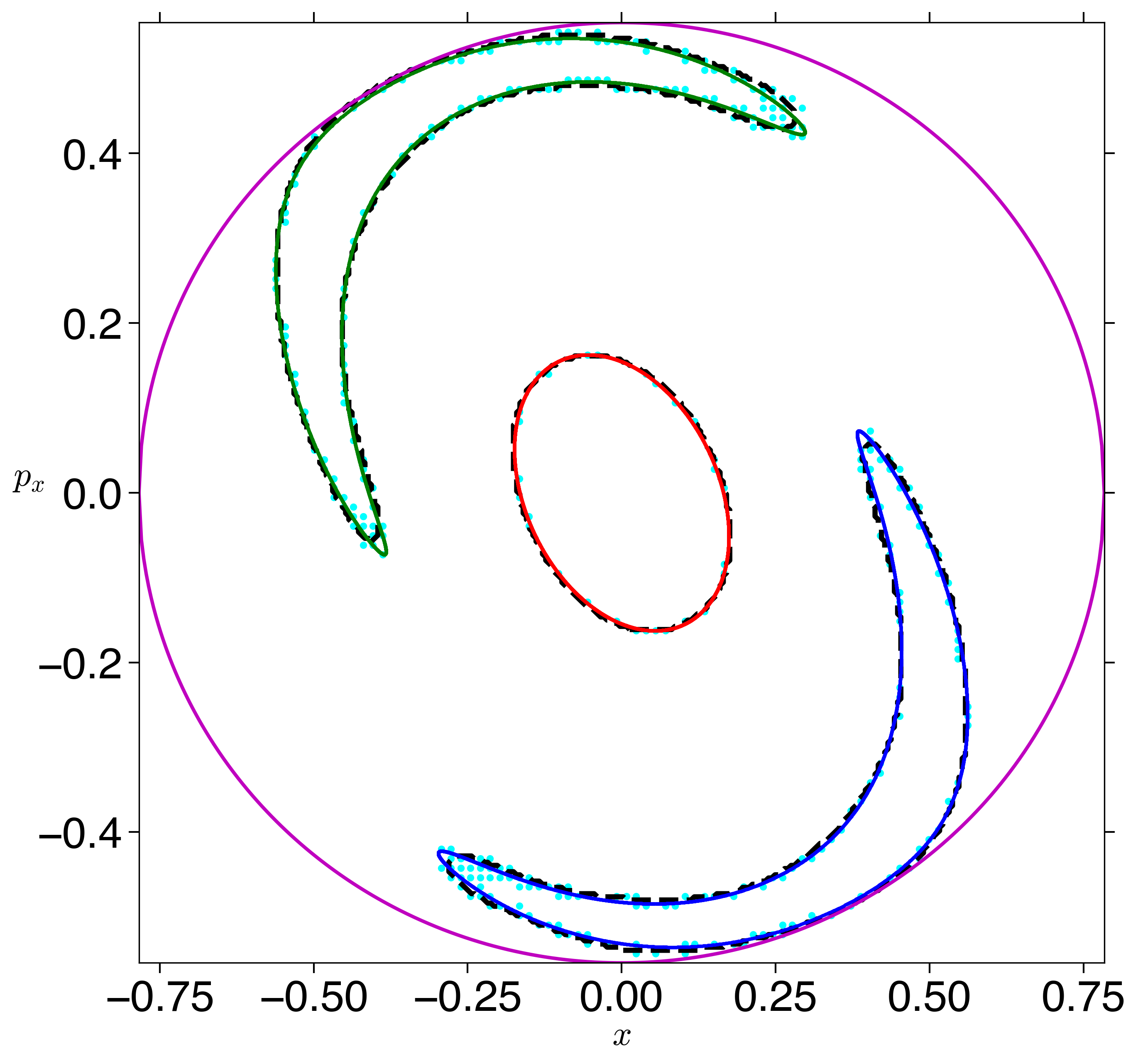}}
	\subfigure[$E = 0.20, y_c = -0.25$]{\includegraphics[width = 0.24\linewidth]{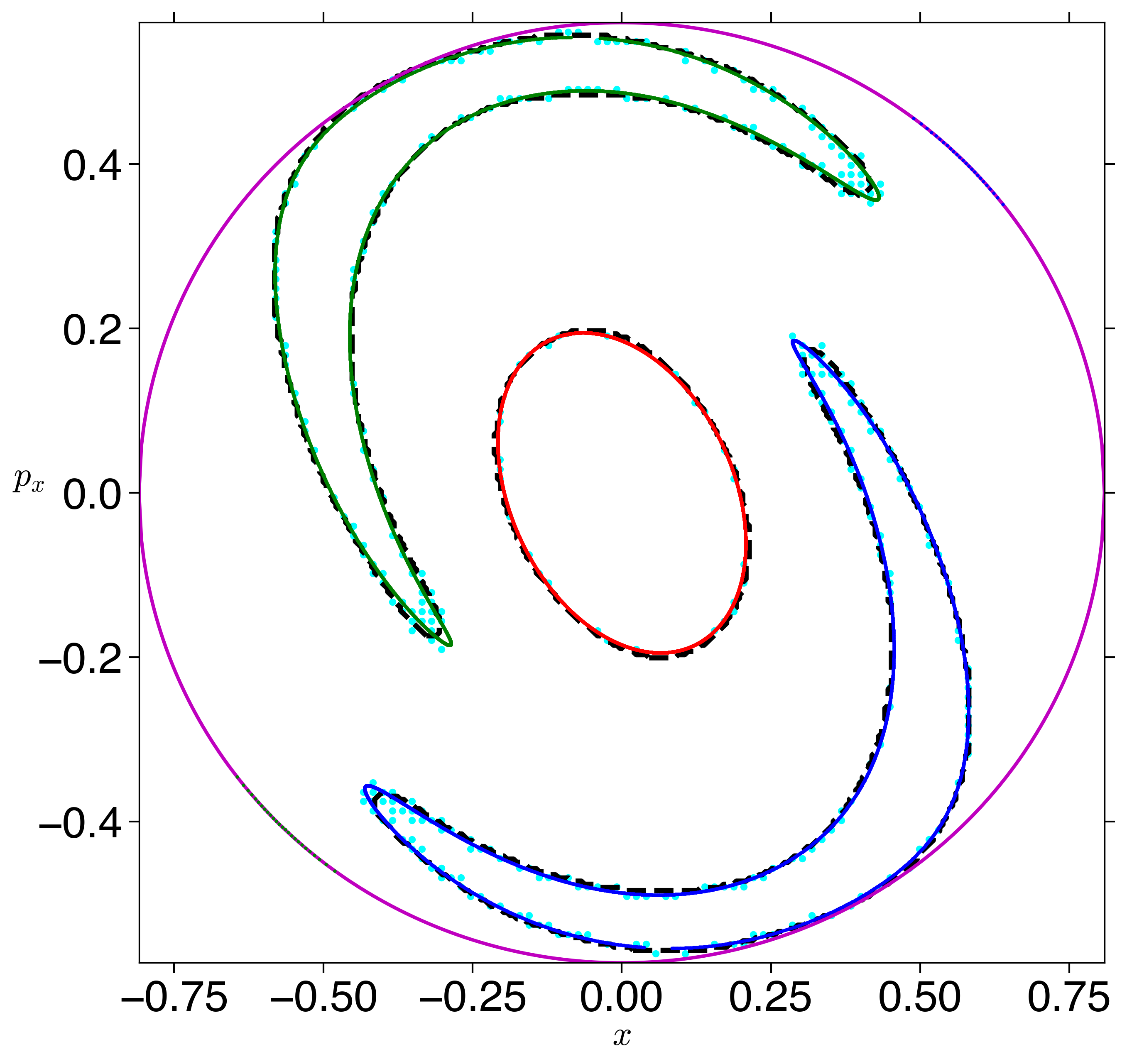}}
	\caption{\textbf{Fixed training data set.} Reactive islands identified by the support vector classifier trained using fixed size data set shown as dashed curves. The overlayed continuous curve is obtained using direct computation of tube manifolds at energy $E = 0.17, 0.18, 0.19, 0.20$ in first, second, third, fourth column, respectively. The magenta curve denotes the intersection of the energy surface with the two dimensional section. The cyan dots denote the support vectors used by the classifier in learning the reactive islands as decision surfaces. Two sections with $(x,p_x)$ coordinates are shown in top and bottom rows: (a-d) $y_c = 0$ (e-h) $y_c = -0.25$ with $p_y > 0$.}
	\label{fig:svc_ris_fixed_td}
\end{figure}
The decision boundaries as learned by the SVC trained using fixed size training dataset and the reactive islands obtained from the direct computation of the tube manifolds are compared for verification in Fig.~\ref{fig:svc_ris_fixed_td}. The learned decision boundaries track the true reactive islands to a high accuracy and the classification accuracy is above 99\% for all the energies considered here. In fact, for the low energy case ($E = 0.17$, which is $0.003$ above the energy of the saddle, the accuracy is similar to the highest energy case, $0.033$. These two energies represent two extreme training data sets as the fraction of reactive trajectories increases with total energy. Thus, for low energy case a uniform sampling is bound to give small number of reactive trajectories and vice versa for high energy case. Even though this can be corrected using a non-uniform sampling for low energy, we show that learning the reactive islands leads to high accuracy because it is the fundamental phase space structure underlying the reactive trajectories.

\textbf{Active learning.} In this approach, we developed an iterative method for {\em active learning} of the reactive islands. The steps in this iterative method are: (i) Begin training step with a relatively small training data that represents a coarse sampling of the two dimensional section using a regular grid. (ii) This training data is then used to train a support vector classifier which learns a coarse decision boundary. (iii) The method adds new training data around the support vectors used by the classifier in the second step. The new data consists of 10 points generated around each support vector using a normal distribution and is added to the initial training data. (iv) The new data set is used in training the classifier which then relearns reactive islands. Steps (iii) and (iv) are repeated until a desired accuracy, for example 99\%, is achieved. For this approach, we use, $C \in \{10, 1e2, 1e3\}$ and $\gamma \in \{1, 10, 1e2\}$, to perform a grid search for highest accuracy along with 5-fold cross-validation during the training. 

\begin{figure}[!h]
	\centering
	\subfigure[$E = 0.17, y_c = 0$]{\includegraphics[width = 0.24\linewidth]{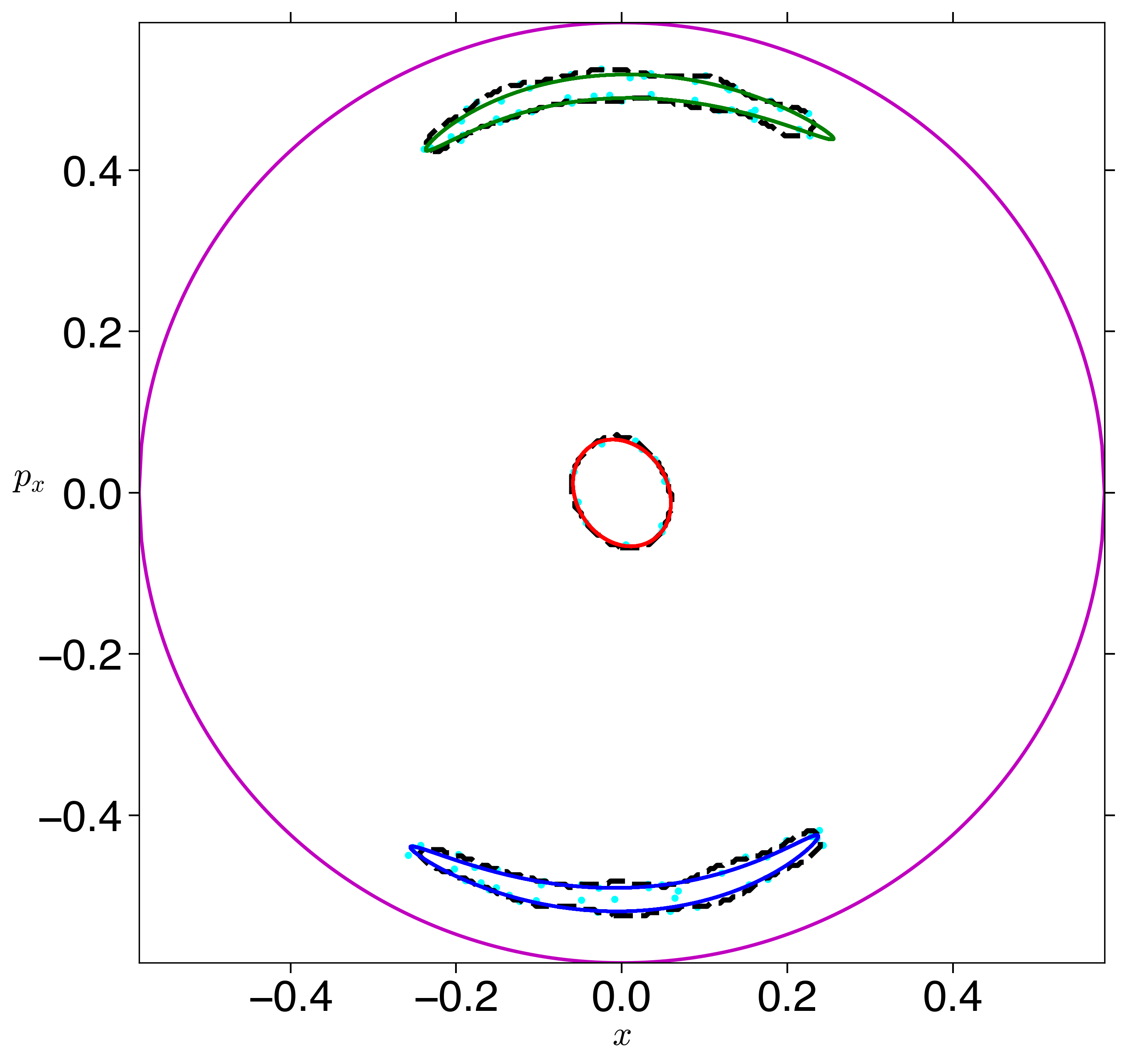}}
	\subfigure[$E = 0.18, y_c = 0$]{\includegraphics[width = 0.24\linewidth]{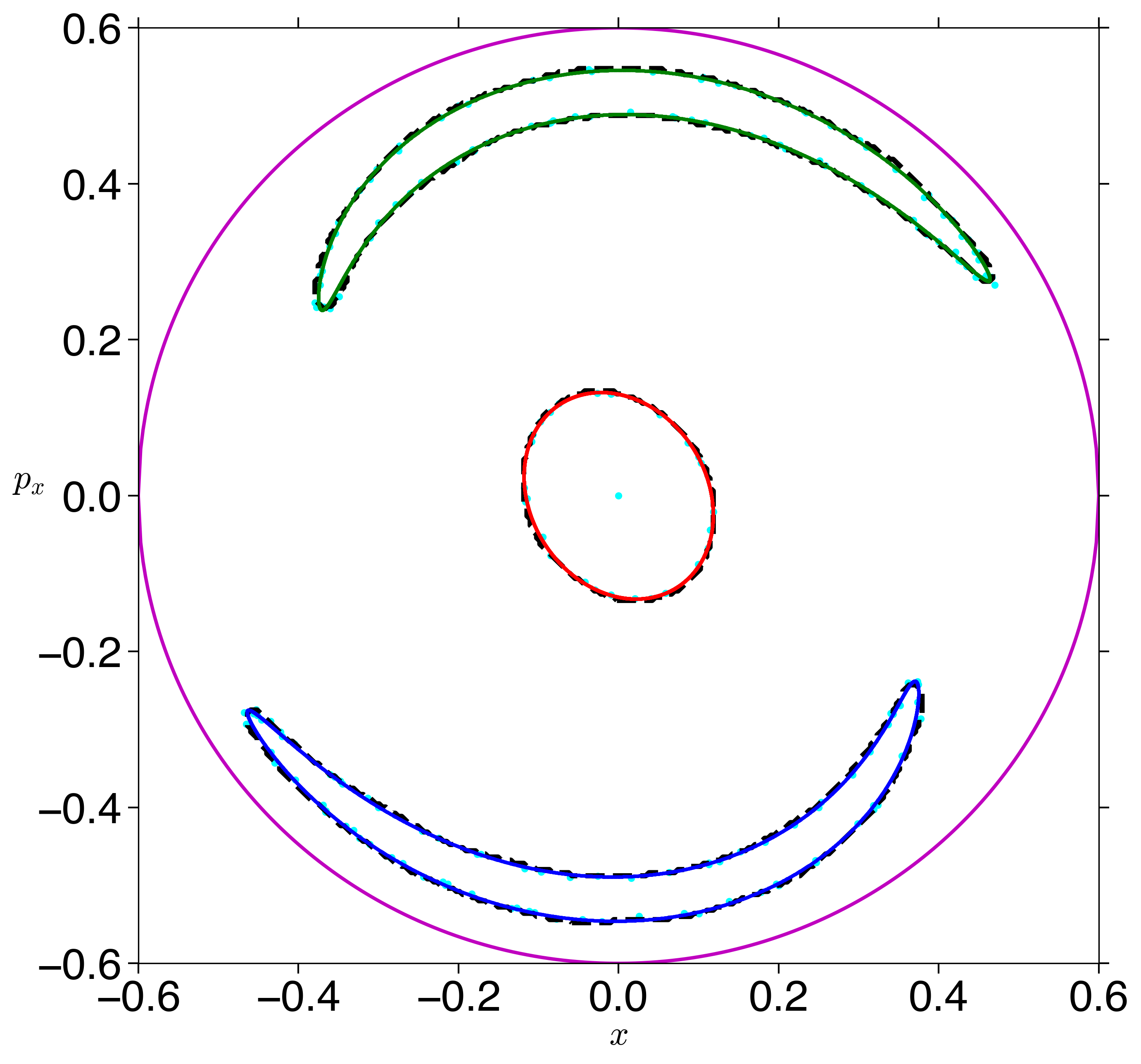}}
	\subfigure[$E = 0.19, y_c = 0$]{\includegraphics[width = 0.24\linewidth]{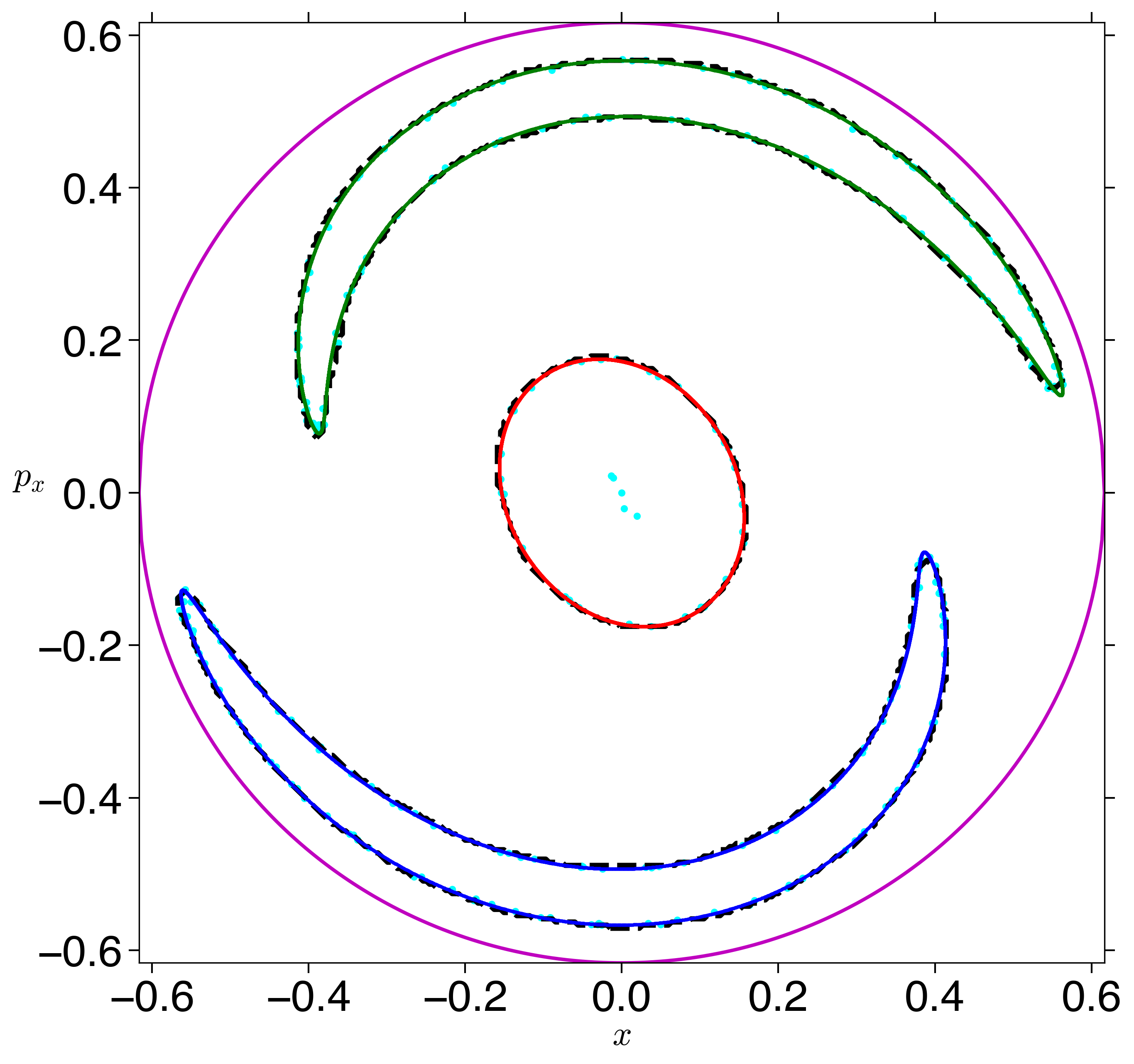}}
	\subfigure[$E = 0.20, y_c = 0$]{\includegraphics[width = 0.24\linewidth]{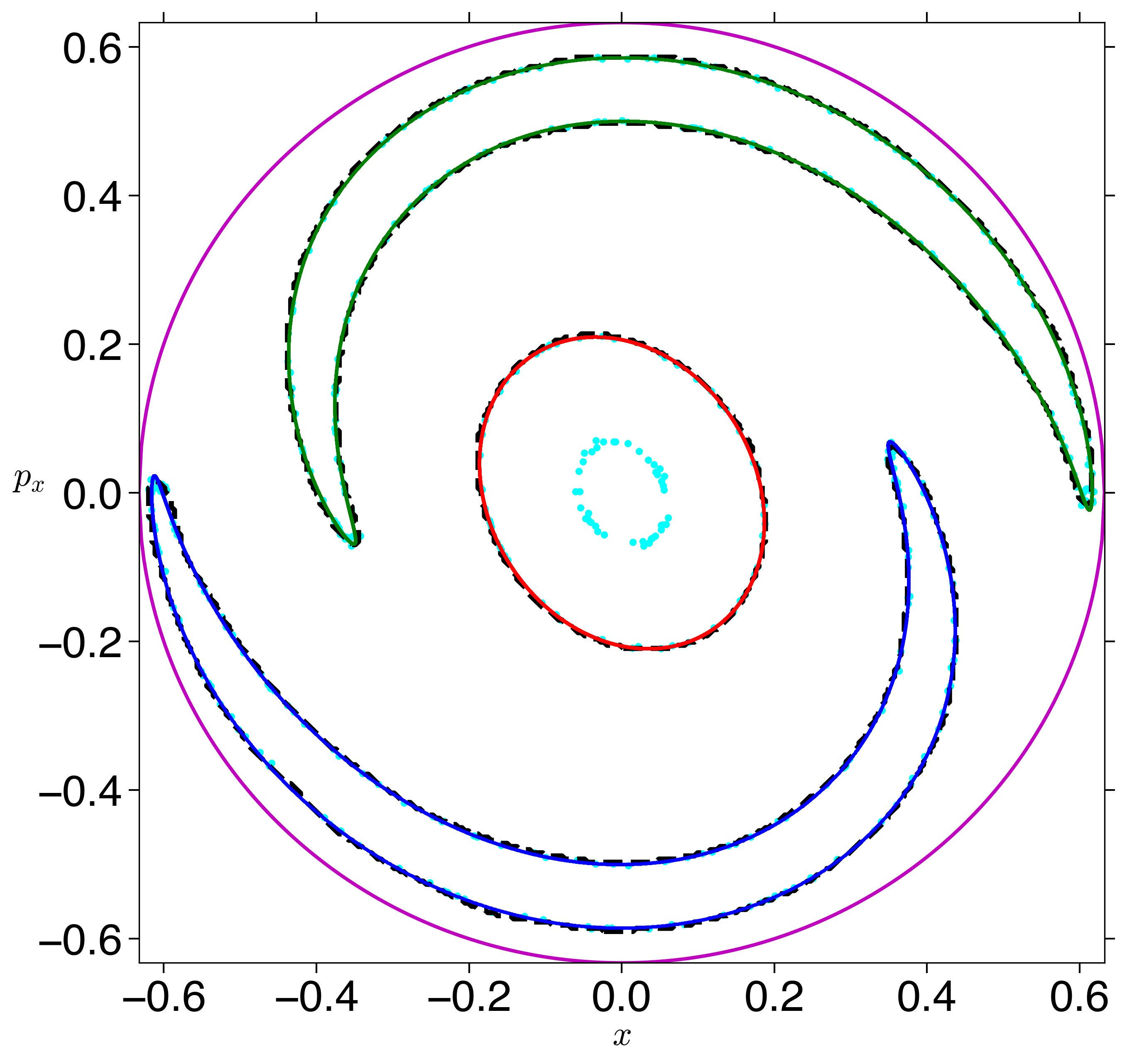}}
	\subfigure[$E = 0.17, y_c = -0.25$]{\includegraphics[width = 0.24\linewidth]{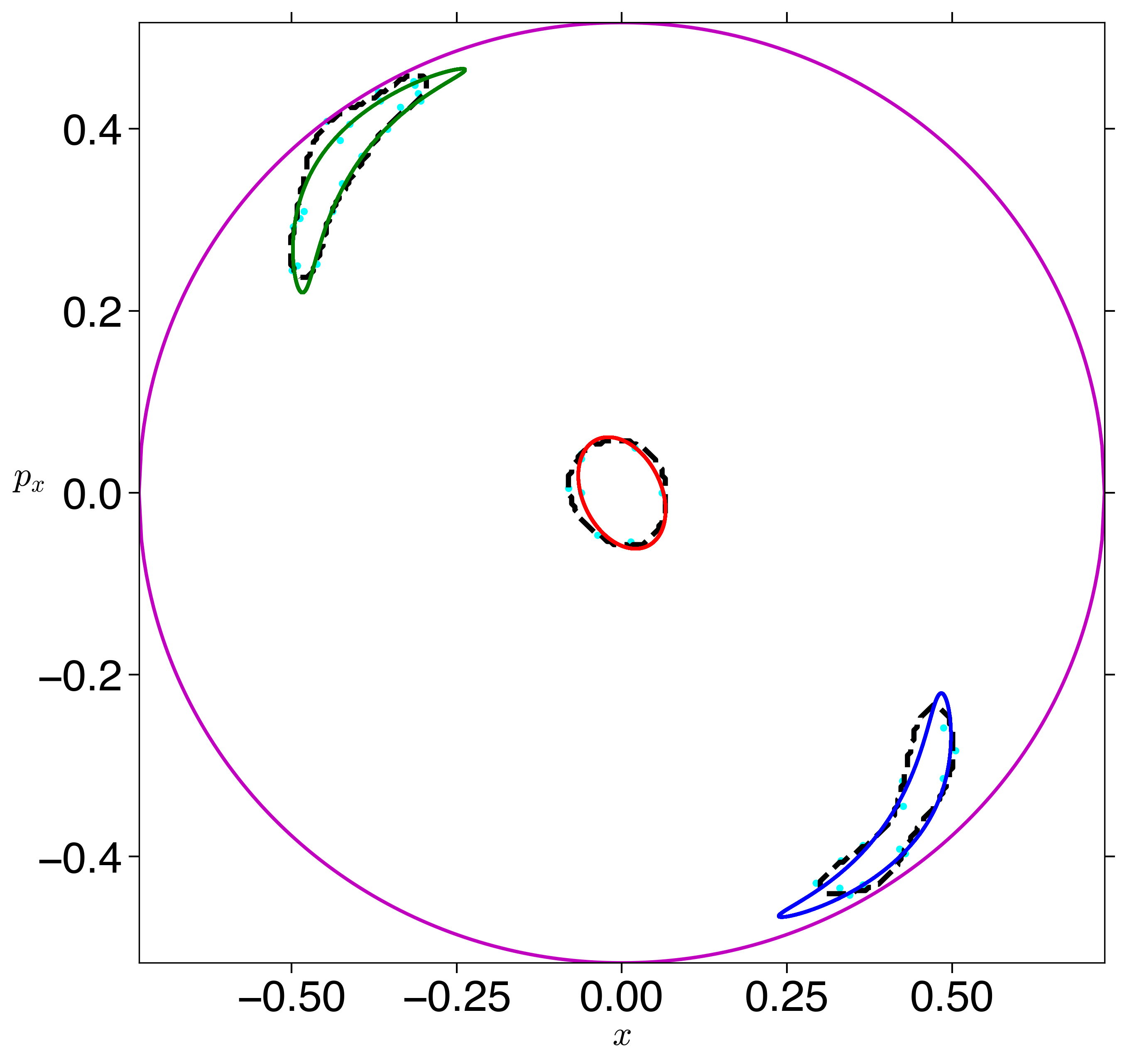}}
	\subfigure[$E = 0.18, y_c = -0.25$]{\includegraphics[width = 0.24\linewidth]{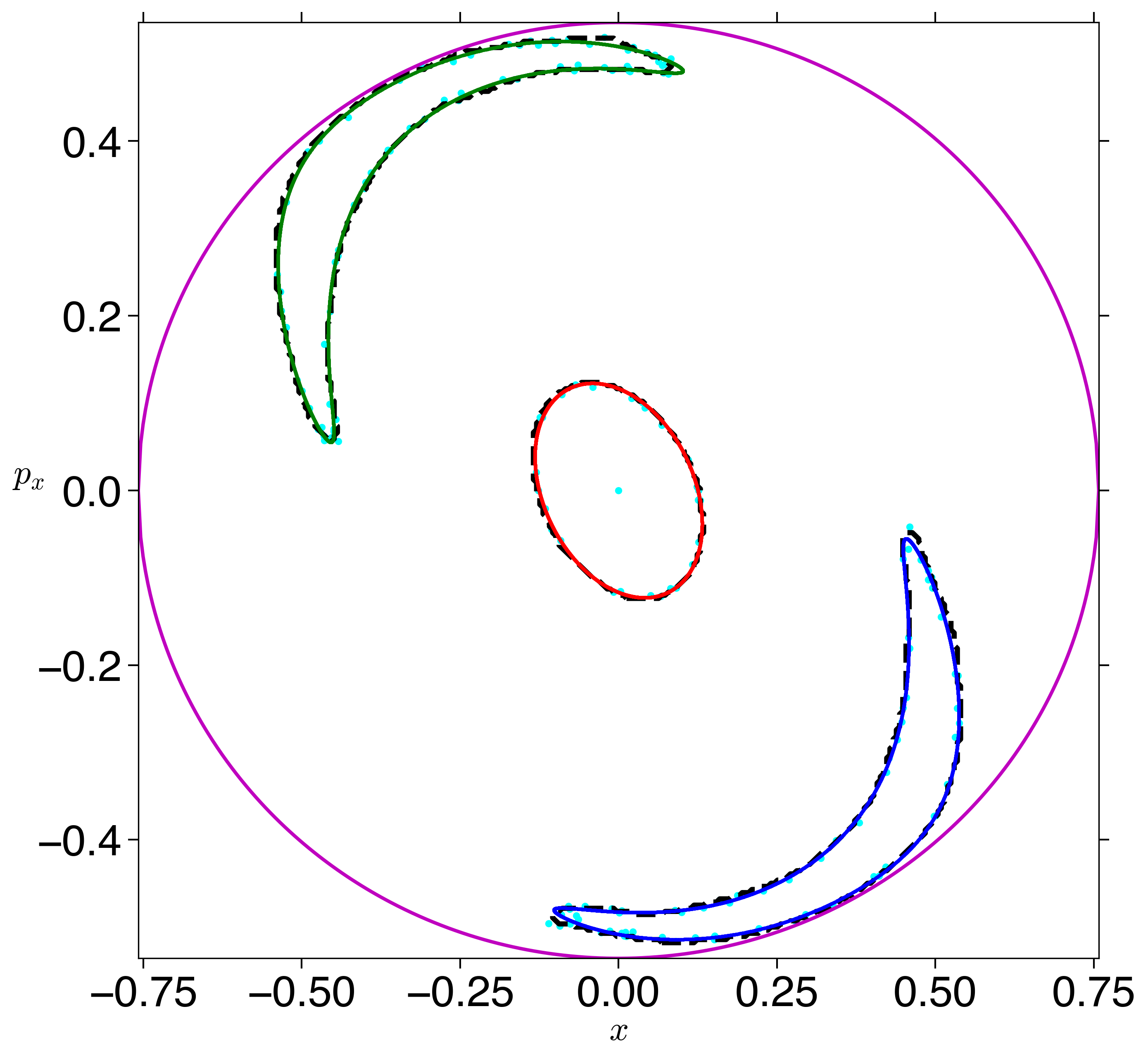}}
	\subfigure[$E = 0.19, y_c = -0.25$]{\includegraphics[width = 0.24\linewidth]{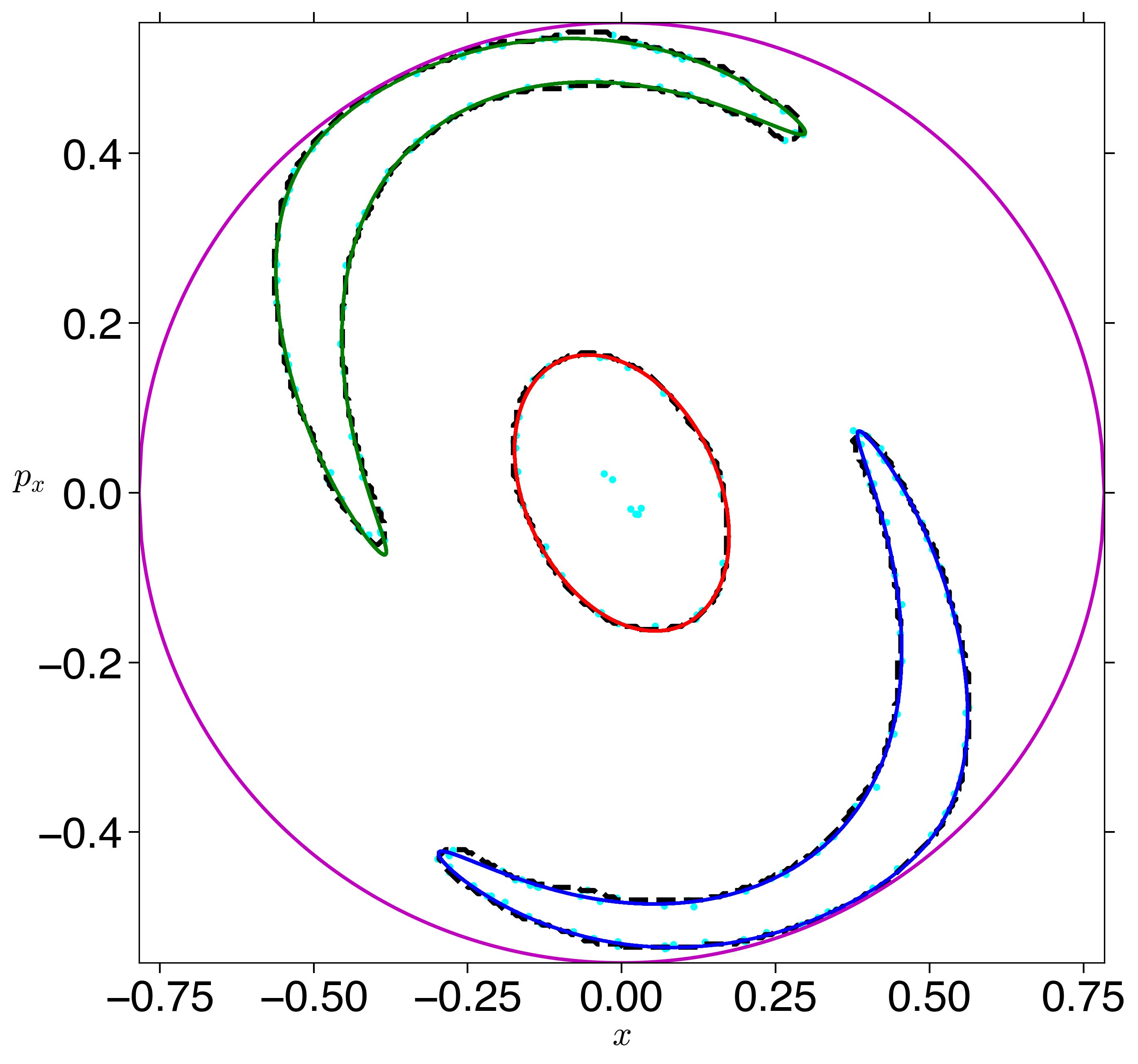}}
	\subfigure[$E = 0.20, y_c = -0.25$]{\includegraphics[width = 0.24\linewidth]{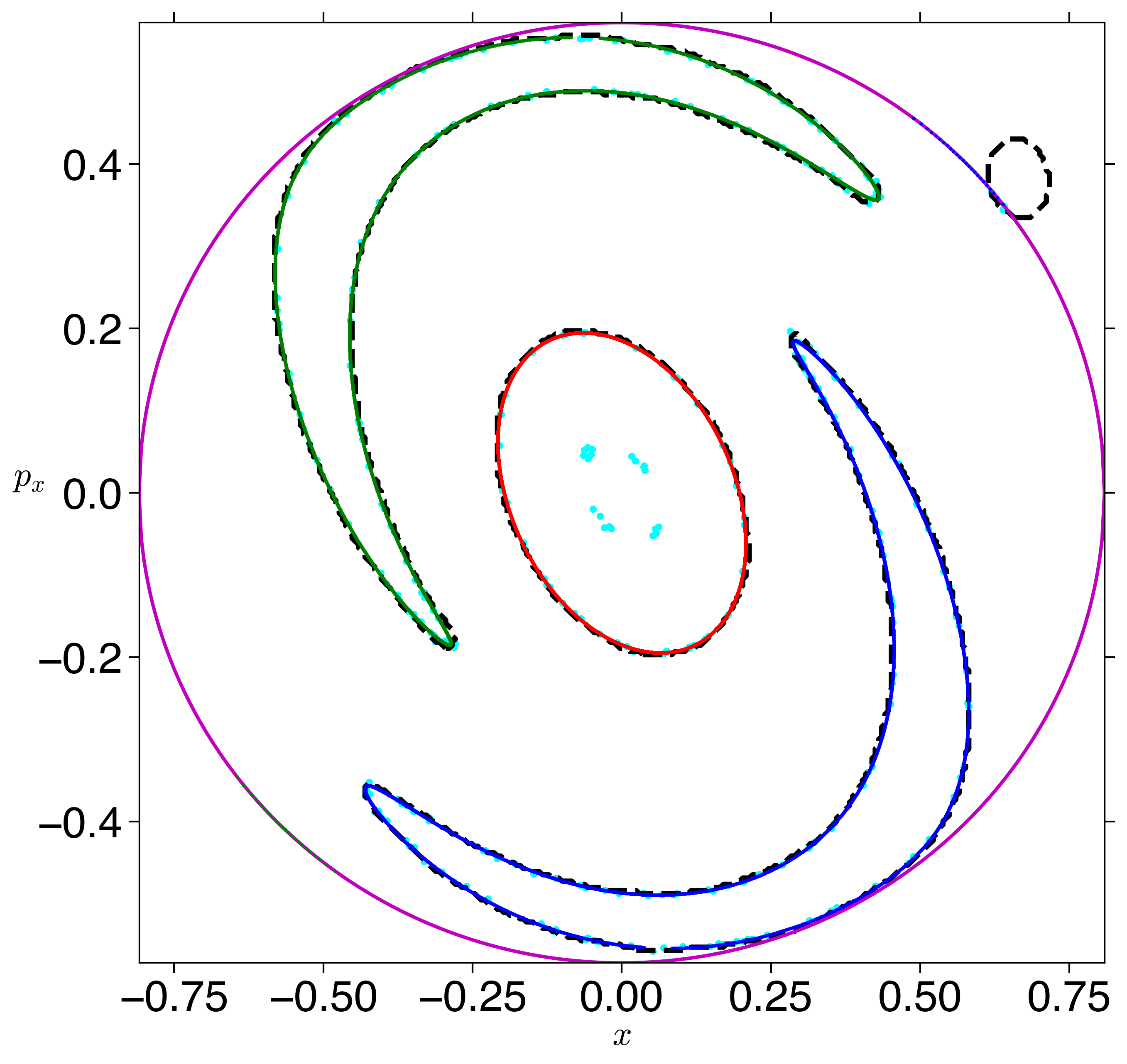}}
	\caption{\textbf{Active learning} Reactive islands identified by the support vector classifier trained using data generated near the coarse boundaries shown as dashed curves. The overlayed continuous curve is obtained using direct computation of tube manifolds at energy $E = 0.17, 0.18, 0.19, 0.20$ in first, second, third, fourth column, respectively. The magenta curve denotes the intersection of the energy surface with the two dimensional section. The cyan dots denote the support vectors used by the classifier in learning the reactive islands as decision surfaces. Two sections with $(x,p_x)$ coordinates are shown in top and bottom rows: (a-d) $y_c = 0$ (e-h) $y_c = -0.25$ with $p_y > 0$.}
	\label{fig:svc_ris_adaptive_td}
\end{figure}

The decision boundaries as learned by the SVC using an active learning approach and the reactive islands obtained from the direct computation of the tube manifolds are compared for verification in Fig.~\ref{fig:svc_ris_adaptive_td}.

\textbf{Trajectory geometry enabled learning.} In this approach, we use a positive scalar quantity for encoding the geometry of a trajectory called the Lagrangian descriptor (LD)~\cite{madrid2009ld,mancho_lagrangian_2013} to generate a new feature for training a SVC. Lagrangian descriptors have been shown to detect phase space structures that mediate transition dynamics in general non-linear dynamical systems (see Ref.~\cite{ld_book_2020} and references therein). Furthermore, as a trajectory based diagnostic of the phase space, it can be computed on-the-fly with the trajectory. This gives it a two-fold merit as a feature: (i) encodes the geometry of the trajectory, thus incorporates the phase space perspective and (ii) efficient computation along with trajectory generation.

We briefly describe the method of Lagrangian descriptors which reveals regions with qualitatively distinct dynamical behavior by showing the intersection of the invariant manifolds with the two dimensional section. For a general time-dependent dynamical system given by
\begin{equation}
	\dfrac{d\mathbf{x}}{dt} = \mathbf{f}(\mathbf{x},t) \;,\quad \mathbf{x} \in \mathbb{R}^{n} \;,\; t \in \mathbb{R} \;,
	\label{eq:gtp_dynSys}
\end{equation}
where the vector field $\mathbf{f}(\mathbf{x},t)$ is assumed to be sufficiently smooth both in space and time. The vector field $\mathbf{f}$ can be prescribed by an analytical model or given from numerical simulations as a discrete spatio-temporal data set. For instance, the vector field could represent the velocity field of oceanic or atmospheric currents obtained from satellite measurements or from the numerical solution of geophysical models. For any initial condition $\mathbf{x}(t_0) = \mathbf{x}_0$, the system of first order nonlinear differential equations given in Eqn.~\eqref{eq:gtp_dynSys} has a unique solution represented by the trajectory that starts from that initial point $\mathbf{x}_0$ at time $t_0$.

\begin{figure}
    \centering
    \subfigure[$E = 0.17, y_c = 0$]{\includegraphics[width = 0.3\linewidth]{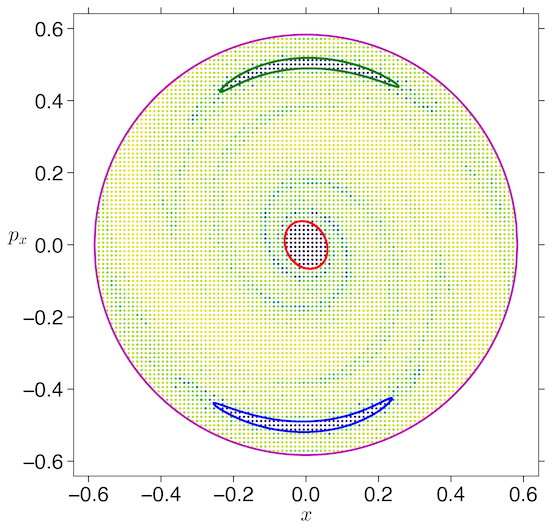}}
    \subfigure[$E = 0.19, y_c = 0$]{\includegraphics[width = 0.3\linewidth]{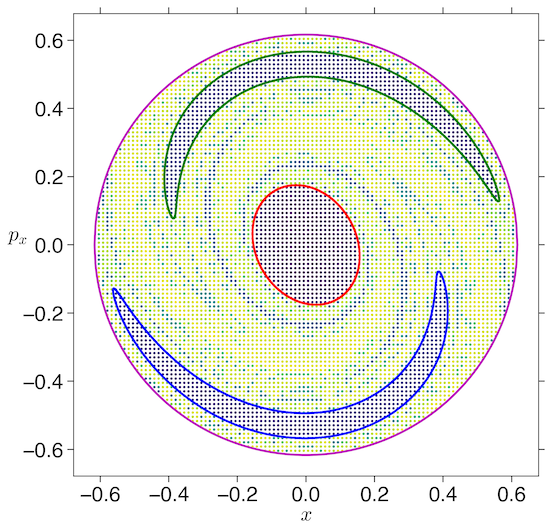}}
	\subfigure[$E = 0.20, y_c = 0$]{\includegraphics[width = 0.3\linewidth]{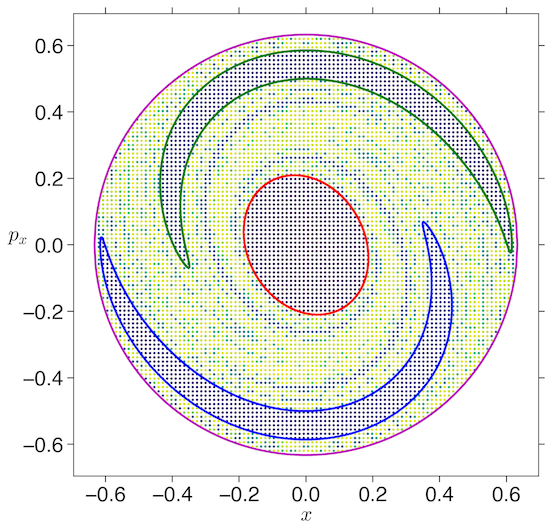}}
    \caption{Reactive islands at $E = 0.18, 0.20$ are shown as red, green, and blue curves, respectively, and the LD values corresponding to the initial conditions sampled on $y_c = 0$ with $p_y > 0$ section. These initial conditions along with the LD value computed for $\tau = 30$ or until a trajectory escapes is used as the training data.}
    \label{fig:ld_ris_vartau30}
\end{figure}

In this study, we adopt the LD definition
\begin{equation}
	\mathcal{L}_p(\mathbf{x}_{0},t_0,\tau) = \int^{t_0+\tau}_{t_0-\tau} \, \sum_{k=1}^{n}   \vert f_{k}(\mathbf{x}(t;\mathbf{x}_0),t) \vert^p \; dt  \;, \quad p \in (0,1]
	\label{eq:Mp_function}
\end{equation}
where $f_k$ is the $k-$the component of the vector field, Eqn.~\eqref{eq:gtp_dynSys} and use $p = 1/2$. We note that the integral can be split into its forward and backward time parts to detect the intersection of stable and unstable manifolds separately. This relates to finding the escape and entry channels into the potential well. In this study, we keep the forward part of the integral given by
\begin{equation}
	\begin{split}
		\mathcal{L}_p^{f}(\mathbf{x}_{0},t_0,\tau) & = \int^{t_0+\tau}_{t_0} \sum_{k=1}^{n} |f_{k}(\mathbf{x}(t;\mathbf{x}_0),t)|^p \; dt
	\end{split}
\end{equation}
Although this definition of LD does not have an intuitive physical interpretation as that of the arclength definition~\cite{mancho_lagrangian_2013}, it allows for a rigorous proof that the ``singular features" (non-differentiable points) in the LD contour map identify intersections with stable and unstable invariant manifolds~\cite{lopesino2017}. Another important aspect of what is known in LD literature as the
$p$-(quasi)norm is that degrees of freedom with relevance in escape/transition (reaction) dynamics can be decomposed and computed. This definition was used to show that the method can be used to successfully detect NHIMs and their stable and unstable manifolds in H{\'e}non-Heiles Hamiltonian~\cite{demian2017,naik2019b}. For this system, where both fixed (or variable) integration time is used, it has also been shown that the LD scalar field attains a minimum (or maximum) value along with singularity at the intersections of the stable and unstable manifolds, and given by
\begin{equation}
	\mathcal{W}^s(\mathbf{x}_{0},t_0) = \textrm{argmin } \mathcal{L}_p^{f}(\mathbf{x}_{0},t_0,\tau) \;,
	\label{eq:min_LD_manifolds}
\end{equation}
where $\mathcal{W}^s(\mathbf{x}_{0},t_0)$ are the stable manifolds calculated at time $t_0$ and  $\textrm{argmin}$ denotes the phase space coordinates on the two dimensional section that minimize the scalar field, $\mathcal{L}_p^{f}(\mathbf{x}_{0},t_0,\tau)$, over the integration time, $\tau$. Thus, the scalar field plotted as a contour map identifies the intersection of the stable manifold with a two dimensional section. This ability of LD contour map to partition trajectories with different phase space geometry is shown in Fig.~\ref{fig:ld_ris_vartau30}(a-c) as values of LD inside the reactive islands are close to constant. 

\begin{figure}[!ht]
	\centering
    \subfigure[$E = 0.17, y_c = 0$]{\includegraphics[width = 0.3\textwidth]{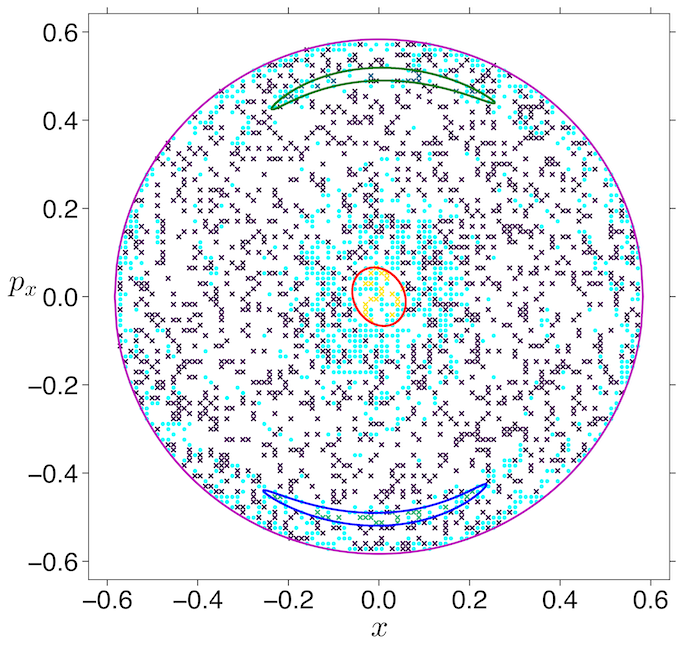}}
	\subfigure[$E = 0.19, y_c = 0$]{\includegraphics[width = 0.3\textwidth]{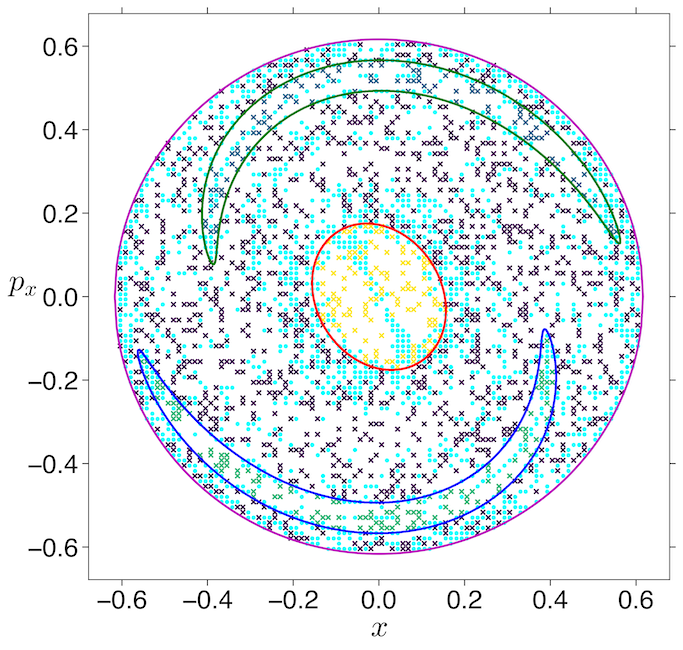}}
	\subfigure[$E = 0.20, y_c = 0$]{\includegraphics[width = 0.3\textwidth]{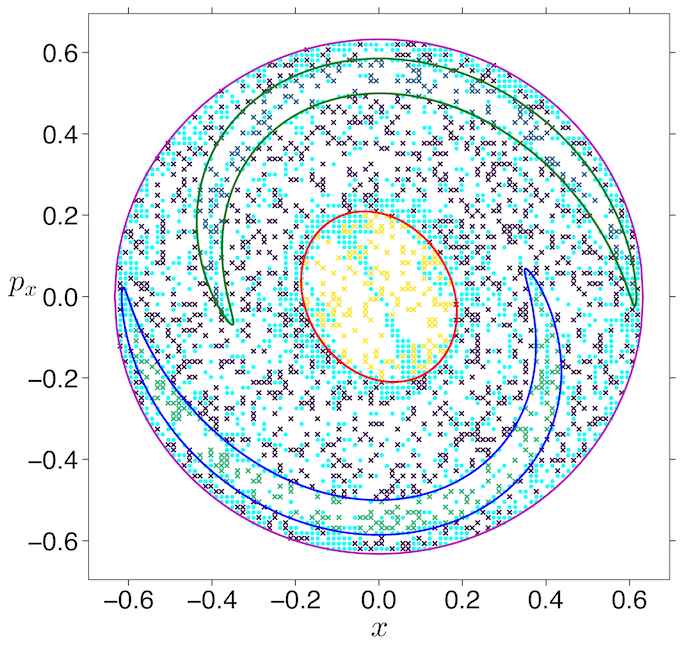}}
	\subfigure[$E = 0.17, y_c = -0.25$]{\includegraphics[width = 0.3\textwidth]{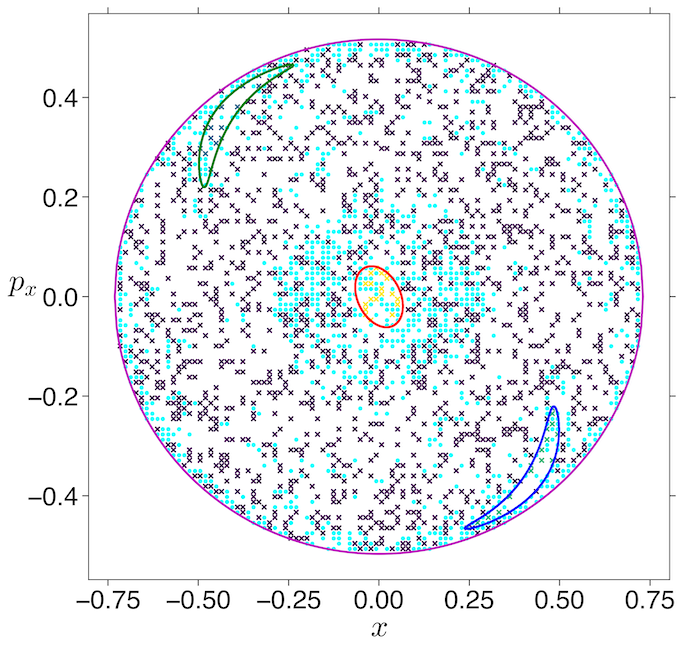}}
	\subfigure[$E = 0.19, y_c = -0.25$]{\includegraphics[width = 0.3\linewidth]{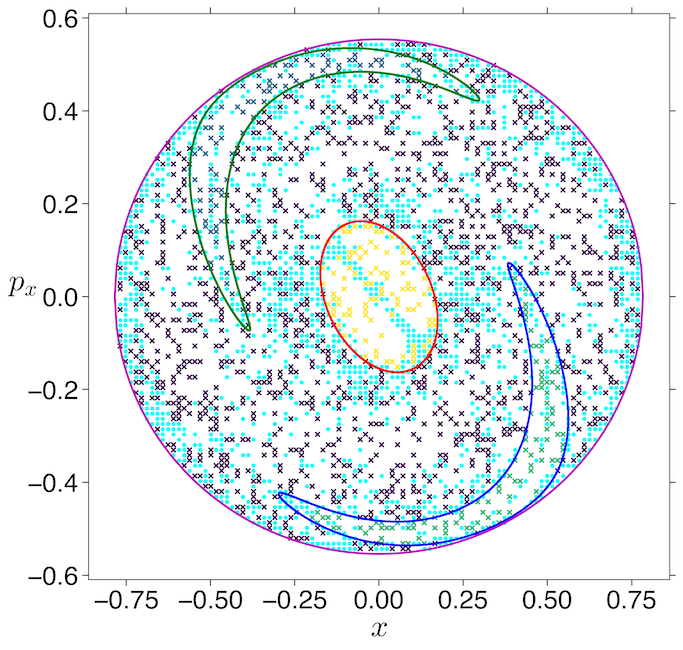}}
	\subfigure[$E = 0.20, y_c= -0.25$]{\includegraphics[width = 0.3\linewidth]{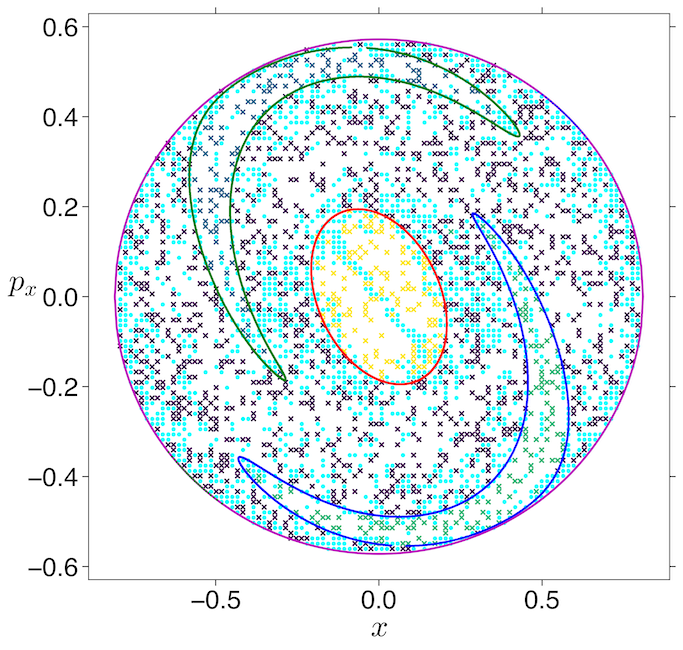}}
	\caption{\textbf{Trajectory geometry enabled learning} Shows the reactive islands, support vectors (as dots), and predictions (as cross) of the trajectory geometry trained support vector classifier.}
	\label{fig:svc_ris_ld_td}
\end{figure}

We construct the training data with three features -  $x,p_x,\mathcal{M}_{0.5}(\tau)$ - and a fixed size dataset. Then, we implement a SVC as in the \emph{Fixed training data} approach with the parameters for the Gaussian radial basis function $C \in \{1e2, 1e3, 1e4, 1e5\}$ and $\gamma \in \{10, 1e2, 1e3, 1e4\}$. In Fig.~\ref{fig:svc_ris_ld_td} (d-i), we show the reactive islands along with the predictions of a SVC trained using the trajectory geometry given by LD as a feature and with a training data size of 10000 points. We note that when such a data set is used the support vectors used by the model increases around the boundary. However, this approach encodes the geometry of a trajectory in phase space and hence leads to robust classification as the total energy parameter is increased.


\section{Conclusions and outlook}
\label{sec:concl}

On a transverse two dimensional section of the energy surface, reactive islands are the intersections of stable and unstable invariant manifolds of hyperbolic periodic orbits. Thus, the one dimensional boundaries of the reactive islands separate transition and non-transition trajectories. In this article, we presented three support vector classifier approaches for learning the reactive islands: fixed dataset training, active learning, and trajectory geometry enabled training. The advantages of our approach are as follows: (a) avoiding the need to compute hyperbolic periodic orbits and the associated invariant manifolds, in favour of finding the reactive islands directly on a surface of section as the boundary between classes of qualitatively different dynamics, (b) minimising computational cost of trajectory calculations by sampling the section near a boundary learned from a coarser sampling (c) using trajectory geometry as a dynamical (phase space) feature in the training data and compressing the high dimensional trajectory into a smaller feature set. Inheriting low data requirements from SVM, our approach is designed to work well for systems where integrating trajectories is expensive and is expected to generalise well systems with more than two degrees of freedom.

Our work intends to simplify the process of finding reactive islands, making it easier to generalise and more accessible to a wider scientific audience. Future work in this direction will involve the application to a model of chemical reaction and examples of high dimensional phase space of a system-bath model.

\section{Acknowledgments}

We acknowledge the support of EPSRC Grant No. EP/P021123/1 and Office of Naval Research (Grant No. N00014-01-1-0769).

\section{Data Availability}

The code that support the findings of this study are openly available on the corresponding author's  \href{https://github.com/Shibabrat/learning-reactive-islands}{GitHub repository}. The data used in training the machine learning model is available from the corresponding author upon reasonable request.

\newpage

\bibliography{SNreac.bib,ri_3dof.bib}

\appendix




\section{Hyperbolic periodic orbit and invariant manifolds}\label{app:tube_manifolds}

\subsection{The Linearized Hamiltonian System} 
\label{ssect:linear}

The linearized equation of motion around the index-one saddle equilibria with 
coordinates $(x_e,y_e,0,0)$ is given by expanding the terms of the Hamiltonian~\eqref{eqn:henon-heiles} and keeping the quadratic terms. After making a coordinate to $(x_e, y_e, 0, 0)$ as the origin, the quadratic Hamiltonian function gives the linear system at the equilibrium point 

%
\begin{equation}
	\begin{pmatrix}
		\dot{x} \\
		\dot{y} \\
		\dot{p_x}\\
		\dot{p_y}
	\end{pmatrix} =  \begin{pmatrix}
		0 & 0 & 1/m_x & 0 \\
		0 & 0 & 0 & 1/m_y \\
		-\omega_x^2 - 2y_e & -2x_e  & 0 & 0 \\
		-2x_e & -\omega_y^2 + 2\delta y_e & 0 & 0
	\end{pmatrix}
	\begin{pmatrix}
		x \\
		y \\
		p_x\\
		p_y
	\end{pmatrix}
	\label{eqn:linearized_hh2dof}
\end{equation}
%


The linearized dynamics near the index-one saddle equilibrium points extend to the full nonlinear system due to Moser's generalization of Lyapunov's theorem.

Let us assume the eigenvector to be $v = [k_1, k_2, k_3, k_4]^T$, and the eigenvalue problem 
becomes $\mathbb{J}v = \beta v$. This gives the expressions
\begin{align}
	k_3/m_x =& \beta k_1 \label{eqn:row1} \\
	k_4/m_y =& \beta k_2 \label{eqn:row2} \\
	(-\omega_x^2 - 2y_e)k_1 - 2 x_e k_2 =& \beta k_3 \label{eqn:row3}\\
	-2 x_e k_1 + (-\omega_y^2 + 2 \delta y_e) k_2 =& \beta k_4 \label{eqn:row4}
\end{align}
Let $k_1 = 1$, then using Eqns.~\eqref{eqn:row1}and~\eqref{eqn:row2} the eigenvector becomes 
$\bigl[1, k_2, \beta m_x, \beta m_y k_2 \bigr]$. 



%
Thus, the eigenvectors corresponding to , $\beta = \pm \lambda, \pm i \omega$, can be written as
\begin{equation}
	\begin{aligned}
		v_{\pm \lambda} = \begin{bmatrix}
			1, & \dfrac{ 2 x_e}{-\omega_y^2 + 2 \delta y_e - \lambda^2 m_y}, & \pm \lambda m_x, & \pm \lambda m_y \dfrac{ 2 x_e}{-\omega_y^2 + 2 \delta y_e - \lambda^2 m_y}
		\end{bmatrix}\\
		v_{\pm i \omega} = \begin{bmatrix}
			1, & \dfrac{ 2 x_e}{-\omega_y^2 + 2 \delta y_e + \omega^2 m_y}, & \pm i \omega m_x, & \pm i \omega m_y \dfrac{ 2 x_e}{-\omega_y^2 + 2 \delta y_e + \omega^2 m_y}
		\end{bmatrix},
	\end{aligned}
\end{equation}
respectively. Thus, the general solution of the linear system near 
the saddle equilibrium point is given by
\begin{equation}
	\begin{aligned}
		\mathbf{x}(t) = & \left\{  x(t), y(t), v_x(t), v_y(t)  \right\} 
		= A_1 e^{\lambda t} v_{\lambda} + A_2 
		e^{-\lambda t} v_{-\lambda} + 2{\rm Re}\left( B e^{i \omega t} v_{i \omega} \right)
		\label{eqn:gen_sol_ind1saddle_2dof}
	\end{aligned}
\end{equation}
with $A_1, A_2$ being real and $B = B_1 + i B_2$ being complex.

\subsection{Computing the hyperbolic periodic orbit and associated tube manifolds at the index-one saddle}
\label{ssect:tube_mani} 
For discussing the geometry, we call the equilibrium with positive y-coordinate $\mathbf{x}_{\rm eq, top}$ and negative y-coordinate $\mathbf{x}_{\rm eq, left}$ and $\mathbf{x}_{\rm eq, right}$. 

\textbf{Select appropriate energy above the critical value~\textemdash~} For computation of manifolds that act as \emph{boundary} between the escape and non-escape trajectories, we select the total energy, $E$, above the energy of the index-one saddle, $E_s$, and thus the excess energy 
$\Delta E = E - E_s > 0$.

\textbf{Differential correction and numerical continuation~\textemdash~} 
We present a procedure which computes the hyperbolic periodic orbits associated with an index-one saddle using a small ``seed'' 
initial conditions obtained from the linearized equations of motion at the index-one saddles. Then, corrects the guess using the linearization about the trajectory obtained by evolving the initial guess in the full nonlinear equations. This is called the differential correction~ (more details can be found in Chapter 4 of the book~\citeauthor{Koon2011} and other application of this method includes celestial mechanics, ship dynamics, models of dissociation and isomerization reaction ~\cite{koon2000heteroclinic,naik2017geometry,naik2019b,naik2020}).


\textit{Guess initial condition of the periodic orbit \textemdash} The linearized equations of motion at an equilibrium  
point gives the initial guess for the differential correction method. Let us select an equilibrium point, 
$\mathbf{x}_{\rm eq, left}$. The linearization yields an eigenvalue problem $Av = \gamma v $, where $A$ is the Jacobian 
matrix in Eqn.~\eqref{eqn:linearized_hh2dof} evaluated at the equilibrium point, $\gamma$ is the eigenvalue, and $v = [k_1, k_2, k_3, 
k_4]^T$ is the corresponding eigenvector.  The idea is to use the complex eigenvalue and the corresponding eigenvector to obtain a guess for the initial condition on the periodic orbit and its period $T_{\rm guess, po}$, which should be close to $2\pi/\omega$ (generalization of Lyapunov's theorem) and increase monotonically with excess energy, $\Delta E$. 


The initial condition for a periodic orbit of x-amplitude, $A_x > 0$ can be computed by letting $A_1 = A_2 = 0$ and $t = 0$ 
in Eqn.~\eqref{eqn:gen_sol_ind1saddle_2dof}, and $B = -A_x/2$ (this choice is made to get rid of factor 2) denotes a small amplitude in the general linear solution. Thus, using the eigenvector in the center projection we can write
\begin{equation}
	\begin{aligned}
		\bar{\mathbf{x}}_{\rm 0,g} =& \begin{pmatrix}
			x_{\rm eq,left},y_{\rm eq,left},0,0
		\end{pmatrix}^T + 
		2Re(B v_{i \omega}) \\
		=& \begin{pmatrix}
			x_{\rm eq,left} - A_x, y_{\rm eq, left} - A_x k_2, 0, 0
		\end{pmatrix}^T
	\end{aligned}
\end{equation}
where we consider, without loss of generality, the left index-one saddle equilibrium point. 

\textit{Differential correction of the initial condition \textemdash}~In this step, we introduce correction to one of the coordinates of the initial guess such that the numerical periodic orbit $\bar{\mathbf{x}}_{\rm po}$ satisfies
\begin{align}
	\left\| \bar{\mathbf{x}}_{\rm po}(T) - 
	\bar{\mathbf{x}}_{\rm po}(0) \right\| < \epsilon
\end{align}
for some tolerance $\epsilon \approx 10^{-6}$. For the H{\'e}non-Heiles Hamiltonian, we hold $x-$coordinate constant and correct the initial guess of the $y-$coordinate, use $p_x-$coordinate for terminating event-based integration, and $p_y-$coordinate to test convergence of the periodic orbit. 
Let us denote the flow map of a differential equation $\mathring{\mathbf{x}} = \mathbf{f}(\mathbf{x})$ with initial condition $\mathbf{x}(t_0) = \mathbf{x}_0$ by $\phi(t;\mathbf{x}_0)$. Thus, the displacement of the final state under a perturbation $\delta t$ is given by 
\begin{align}
	\delta \bar{\mathbf{x}}(t + \delta t) = \phi(t + \delta 
	t;\bar{\mathbf{x}}_0 + \delta \bar{\mathbf{x}}_0) - 
	\phi(t ;\bar{\mathbf{x}}_0)
\end{align}
with respect to the trajectory
$\bar{\mathbf{x}}(t)$. Thus, measuring the displacement 
at $t_1 + \delta t_1$ and expanding into Taylor series 
gives
\begin{align}
	\delta \bar{\mathbf{x}}(t_1 + \delta t_1) = 
	\frac{\partial \phi(t_1;\bar{\mathbf{x}}_0)}{\partial 
		\mathbf{x}_0}\delta \bar{\mathbf{x}}_0 + \frac{\partial 
		\phi(t_1;\bar{\mathbf{x}}_0)}{\partial t_1}\delta t_1 + 
	h.o.t
\end{align}
where the first term on the right hand side is the state 
transition matrix, $\mathbf{\Phi}(t_1,t_0)$, when 
$\delta 
t_1 = 0$ and can be obtained by solving the variational equations numerically along with the trajectory~\cite{Parker1989}. Let 
us suppose we want to reach the desired point 
$\mathbf{x}_{\rm d}$, we have
\begin{align}
	\bar{\mathbf{x}}(t_1) = \phi(t_1;\bar{\mathbf{x}}_0) 
	= \bar{\mathbf{x}}_1 = \mathbf{x}_d - \delta 
	\bar{\mathbf{x}}_1
\end{align}
which has an error $\delta \bar{\mathbf{x}}_1$ and needs 
correction. This correction to the first order can be 
obtained from the state transition matrix at $t_1$ and gives a new guess of the periodic orbit. This new initial condition can then be evolved and corrected as an iterative procedure with ``small'' (first order or differential) correction, this gives convergence in few steps. For the index-one saddle equilibrium points in the H{\'e}non-Heiles Hamiltonian, we initialize the guess as
\begin{align}
	\bar{\mathbf{x}}(0) = (x_{0,g},y_{0,g},0,0)^T
\end{align}
and using numerical integrator we obtain the orbit until next the half-period event $p_x = 0$ a high tolerance 
(typically $10^{-14}$). So, we obtain 
$\bar{\mathbf{x}}(t_1)$ which for the guess periodic 
orbit denotes the half-period point, $t_1 = T_{0,g}/2$ 
and compute the state transition matrix 
$\mathbf{\Phi}(t_1,0)$. Using the entries of the state transition matrix, we derive the correction for the coordinate $y_{0,g}$ assuming $x_{0,g}$ is constant. Thus, the first order correction is given by
\begin{align}
	\delta p_{x_1} = \Phi_{32}\delta y_0 + \dot{p}_{x_1}\delta t_1 + h.o.t \\
	\delta p_{y_1} = \Phi_{42}\delta y_0 + \dot{p}_{y_1}\delta t_1 + h.o.t
\end{align}
where $\Phi_{ij}$ is the $(i,j)^{th}$ entry of $\mathbf{\Phi}(t_1,0)$ and the acceleration terms come from the equations of motion evaluated at the crossing 
$t 
= t_1$ when $p_{x_1} = \delta p_{x_1} = 0$. Thus, we 
obtain the first order correction $\delta y_0$ as
\begin{align}
	\delta y_0 &\approx \left(\Phi_{42} - 
	\Phi_{32}\frac{\dot{p}_{y_1}}{\dot{p}_{x_1}} 
	\right)^{-1} \delta p_{y_1} \\
	y_0 &\rightarrow y_0 - \delta y_0
\end{align}
which is iterated until $|p_{y_1}| = |\delta p_{y_1}| < 
\epsilon$ for some tolerance $\epsilon$, since we want 
the final periodic orbit to be of the form 
\begin{align}
	\bar{\mathbf{x}}_{t_1} = (x_1,y_1,0,0)^T
\end{align}
This procedure yields an accurate initial condition for 
a periodic orbit of small amplitude $A_x \approx 10^{-4}$, since 
our initial guess is based on the linearization at 
the equilibrium point.

\textit{Numerical continuation to periodic orbit at arbitrary energy.\textemdash}~The procedure described above yields an 
accurate initial condition 
for a periodic orbit from a single initial guess. If our initial guess came 
from the linear approximation near the equilibrium point, from
Eqn.~\eqref{eqn:gen_sol_ind1saddle_2dof}, it has been observed numerically that we can only use this procedure for small 
amplitude, of order $10^{-4}$, periodic orbits 
around $\mathbf{x}_{\rm eq, bot}$. This small amplitude correspond to small excess 
energy, typically of the order $10^{-2}$, and if we want to compute the periodic orbit of arbitrarily large amplitude, we 
resort to numerical continuation for generating a family which reaches the appropriate total energy. This is done using two 
nearby periodic orbits of small amplitude to obtain initial guess for the next periodic orbit and performing differential 
correction to this guess. To this end, we proceed as follows. Suppose we find two small nearby periodic orbit initial 
conditions, $\bar{\mathbf{x}}_0^{(1)}$ and $\bar{\mathbf{x}}_0^{(2)}$, 
correct to within the tolerance $d_{\rm tol}$, using the differential correction 
procedure described above. We can generate a family of periodic orbits with 
successively increasing amplitudes around $\bar{\mathbf{x}}_{\rm eq, bot}$ in 
the following way. Let 
\begin{equation}
	\Delta = \bar{\mathbf{x}}_0^{(2)} - \bar{\mathbf{x}}_0^{(1)} 
	= [\Delta x_0, \Delta y_0, 0, 0]^T		
\end{equation}
A linear extrapolation to an initial guess of slightly larger amplitude, 
$\bar{\mathbf{x}}_0^{(3)}$ is given by
\begin{equation}
	\begin{aligned}
		\bar{\mathbf{x}}_{0,g}^{(3)} =&~\bar{\mathbf{x}}_0^{(2)} + \Delta \\
		=& \left[(\mathbf{x}_0^{(2)} + \Delta x_0), (y_0^{(2)} + \Delta y_0), 0, 0  
		\right] ^T \\
		=& \left[x_0^{(3)}, y_0^{(3)}, 0, 0  \right] ^T
	\end{aligned}
\end{equation}
Thus, keeping $x_0^{(3)}$ fixed, we can use differential correction on this 
initial condition to compute an accurate solution $\bar{\mathbf{x}}_0^{(3)}$ 
from the initial guess $\bar{\mathbf{x}}_{\rm 0,g}^{(3)}$ and repeat the 
process until we have a family of solutions. We can keep track of the energy of 
each periodic orbit and when we have two solutions, $\bar{\mathbf{x}}_0^{\rm 
	(k)}$ and $\bar{\mathbf{x}}_0^{\rm (k+1)}$, 
whose energy brackets the appropriate energy, $E$, we can resort to combining 
bisection and differential correction to these two periodic orbits until we 
converge to the desired periodic orbit to within a specified 
tolerance. Thus, the result is a periodic orbit at desired total energy $E$ and of 
some period $T$ with an initial condition $X_0$.

\begin{figure}[!ht]
	\centering
	\includegraphics[width=0.9\textwidth]{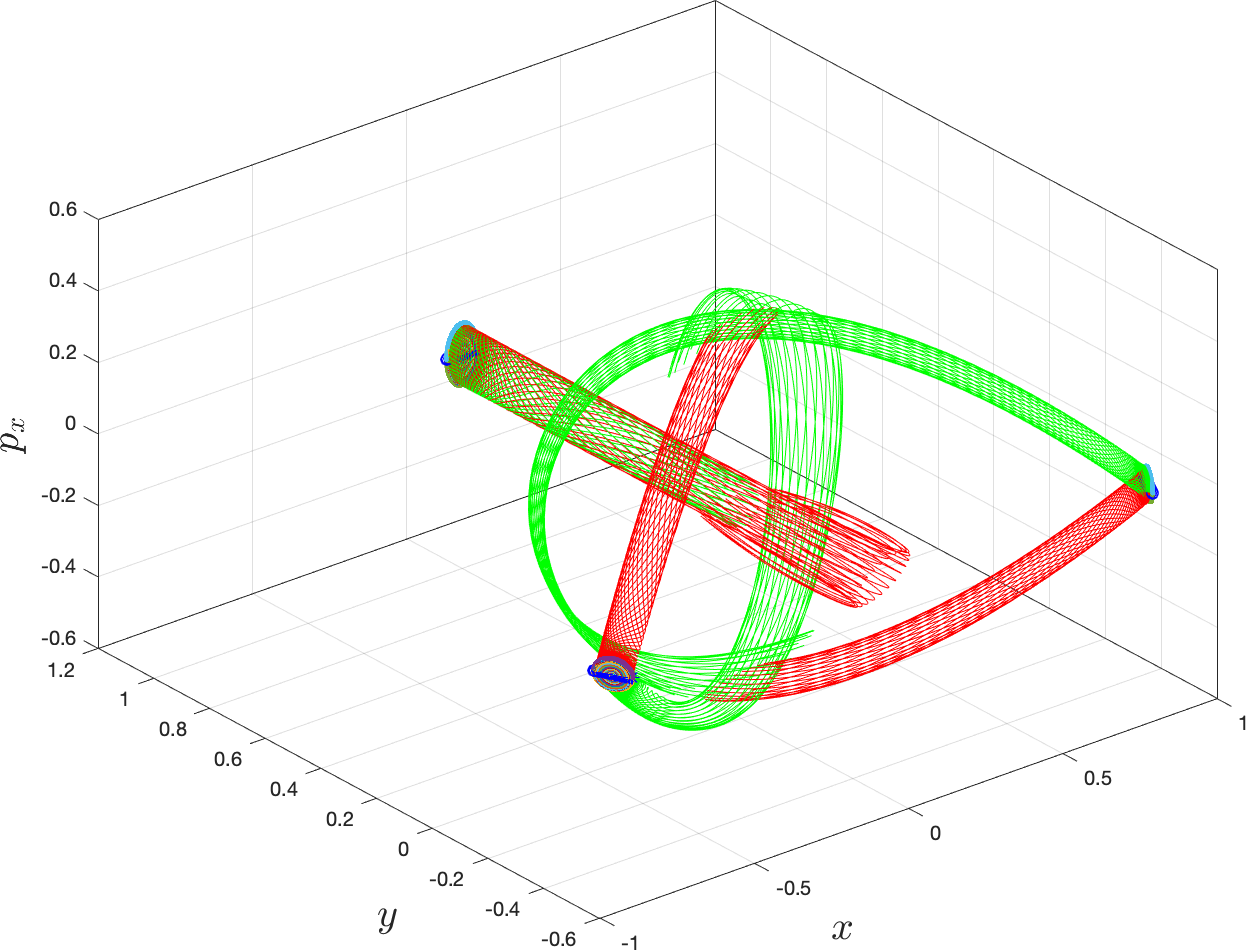}
	\caption{Tube manifolds for the three index-one saddles at total energy $E=0.17$. The green and red trajectories denote the stable and unstable tube manifolds, respectively. We only show the branches of the manifolds inside the potential well. The branches of the stable manifolds mediate escape out of the well and branches of the unstable manifolds mediate entry into the well.}
	\label{fig:manifolds_all_saddles}
\end{figure}

\textbf{Globalization of invariant manifolds \textemdash}~We find the local 
approximation to the 
unstable and stable manifolds of the periodic orbit from the eigenvectors of 
the monodromy matrix. 
Next, the local linear approximation of the unstable (or stable) manifold in 
the form of a state vector is integrated in the nonlinear equations of motion 
to produce the approximation of the unstable (or stable) manifolds, such as those shown in Fig. \ref{fig:manifolds_all_saddles}.
This procedure is known as \textit{globalization of the manifolds} and we proceed as 
follows.

First, the state transition matrix $\Phi(t)$ along the periodic orbit with initial condition $X_0$ can be obtained 
numerically by integrating the variational equations 
along with the equations of motion from $t = 0$ to $t = T$. This is known as the monodromy matrix $M = \Phi(T)$ and the 
eigenvalues can be computed numerically. 
For 
Hamiltonian systems (see~\cite{Meyer2009} for details), tells us that the four 
eigenvalues of $M$ are of the form
\begin{align}
	\lambda_1 > 1, \qquad \lambda_2 = \frac{1}{\lambda_1}, \qquad \lambda_3 = 
	\lambda_4 = 1
\end{align}
The eigenvector associated with eigenvalue $\lambda_1$ is in the unstable 
direction, 
the eigenvector associated with eigenvalue $\lambda_2$ is in the stable 
direction. Let 
$e^{s}(X_0)$ denote the normalized (to 1) stable eigenvector, and $e^{u}(X_0)$ 
denote 
the normalized unstable eigenvector. We can compute the manifold by 
initializing along 
these eigenvectors as:
\begin{equation}
	X^s(X_0) = X_0 + \epsilon e^s(X_0)
\end{equation}
for the stable manifold at $X_0$ along the periodic orbit as
\begin{equation}
	X^u(X_0) = X_0 + \epsilon e^u(X_0)
\end{equation}
for the unstable manifold at $X_0$. Here the small displacement from $X_0$ is 
denoted 
by $\epsilon$ and its magnitude 
should be small enough to be within the validity of the linear estimate, yet 
not so 
small that the time of flight becomes too large due to asymptotic nature of the 
stable 
and unstable manifolds. Ref.~\cite{Koon2011} suggests typical values of 
$\epsilon > 
0$ corresponding to nondimensional position displacements of magnitude around 
$10^{-6}$. 
By numerically integrating the unstable vector forwards in time, using both 
$\epsilon$ 
and $-\epsilon$, for the forward and backward branches respectively, we 
generate 
trajectories shadowing the two branches, $W^u_{+}$ and $W^u_{-}$, of the 
unstable 
manifold of the periodic orbit. Similarly, by integrating the stable vector 
backwards 
in time, using both $\epsilon$ and $-\epsilon$, for forward and backward branch 
respectively, we generate trajectories shadowing the stable manifold, 
$W^{s}_{+,-}$. 
For the manifold at $X(t)$, one can simply use the state transition matrix to 
transport the eigenvectors from $X_0$ to $X(t)$ as
\begin{equation}
	X^s(X(t)) = \Phi(t,0)X^s(X_0)
\end{equation}
It is to be noted that since the state transition matrix does not preserve the 
norm, the resulting vector must be normalized. The globalized invariant 
manifolds associated with index-one saddles are known as 
Conley-McGehee tubes~\cite{Marsden2006}. These tubes form the skeleton 
of transition dynamics by acting as conduits for the states inside them to 
travel between potential wells.
%
%

The computation of codimension-1 separatrix associated with the hyperbolic periodic orbit around a index-one saddle begins with 
the linearized equations of motion. This is obtained after a coordinate transformation to the saddle equilibrium point and 
Taylor expansion of the equations of motion. Keeping the first order terms in this expansion, we obtain the eigenvalues and 
eigenvectors of the linearized system. The eigenvectors corresponding to the center direction provide the starting guess for 
computing the hyperbolic periodic orbits for small excess energy, $\Delta E << 1$, above the saddle's energy. This iterative 
procedure performs small correction to the starting guess based on the terminal condition of the periodic orbit until a 
desired tolerance is satisfied. This procedure is known as differential correction and generates hyperbolic periodic orbits 
for small excess energy. Next, a numerical continuation is implemented to follow the small energy (amplitude) periodic 
orbits out to high excess energies. 

\end{document}